\documentclass[twoside,letterpaper,titlepage,11pt]{book}
\usepackage{graphicx}
\usepackage{amsmath}
\usepackage{amsfonts}
\usepackage{amssymb}
\usepackage{verbatim}
\usepackage{color}
\usepackage{chap}
\usepackage[multiple]{footmisc}

\newcommand{\cD}{{\mathcal D}}
\newcommand{\R}{{\mathbb R}}

\newcommand{\N}{{\mathbb N}}
\newcommand{\cK}{{\mathcal K}}

\newcommand{\myboxed}[1]{\vspace{5 mm}\par \noindent
\framebox[\textwidth][c]{\begin{minipage}[c]{0.95 \textwidth}
#1 \end{minipage}}\vspace{5 mm}\par}

\newcommand{\name}{SOSTOOLS}

\newenvironment{matlab}{\vspace{2 mm}\par
\begin{minipage}[c]{0.90 \textwidth} \tt
}{
\end{minipage}\vspace{3 mm}\par\noindent\hspace{-0.4em}}

\begin{document}
\vspace{-2.5cm}
\title{\textbf{\huge \name}\\ $\,$  \\Sum of Squares Optimization Toolbox for MATLAB\\ $\,$ \\User's guide
\vspace{2em}\\ Version 4.00\vspace{0em}\\ {\normalsize 14th September 2021}}
\date{}
\author{\vspace{0.05cm}\\Antonis Papachristodoulou\footnotemark[1]\qquad James Anderson\footnotemark[2]\qquad Giorgio Valmorbida\footnotemark[3]\\ 
 Stephen Prajna\footnotemark[4] \qquad Peter Seiler\footnotemark[5] \qquad Pablo A. Parrilo\footnotemark[6]\\
\vspace{1cm} Matthew M. Peet\footnotemark[7] \qquad Declan Jagt\footnotemark[7]\vspace{0.15cm}\\
\footnotemark[1]Department of Engineering Science\\University of Oxford, Oxford, U.K.\vspace{0.15cm}\\
\footnotemark[2]Columbia University \\ New York, NY -- USA\vspace{0.15cm} \\
\footnotemark[3]Laboratoire de Signaux et Systèmes \\ CentraleSup{\'e}lec, Gif sur Yvette, 91192, France\vspace{0.15cm} \\
\footnotemark[4]Control and Dynamical Systems\\
California Institute of Technology, Pasadena, CA 91125 -- USA\vspace{0.15cm}\\
\footnotemark[5]Aerospace and Engineering Mechanics Department\\
University of Minnesota, Minneapolis, MN 55455-0153 -- USA\vspace{0.15cm}\\
\footnotemark[6]Laboratory for Information and Decision Systems\\
Massachusetts Institute of Technology, Massachusetts, MA  02139-4307 -- USA\vspace{0.15cm}\\
\footnotemark[7]School for the Engineering of Matter, Transport, and Energy\\
Arizona State University, Tempe, AZ  85255-6106 -- USA\vspace{-.3cm}
}
\maketitle
\noindent Copyright (C) 2002, 2004, 2013, 2016, 2018, 2021  A.\ Papachristodoulou, J.\ Anderson, G.\ Valmorbida,  S.\ Prajna, P.\ Seiler, P.\ A.\ Parrilo, M.\ Peet, D.\ Jagt \vspace{2em}

\noindent This program is free software; you can redistribute it and/or
modify it under the terms of the GNU General Public License
as published by the Free Software Foundation; either version 2
of the License, or (at your option) any later version.\vspace{2em}

\noindent This program is distributed in the hope that it will be useful,
but WITHOUT ANY WARRANTY; without even the implied warranty of
MERCHANTABILITY or FITNESS FOR A PARTICULAR PURPOSE.  See the
GNU General Public License for more details.\vspace{2em}

\noindent You should have received a copy of the GNU General Public License
along with this program; if not, write to the Free Software
Foundation, Inc., 59 Temple Place - Suite 330, Boston, MA  02111-1307, USA.\vspace{2em}

\tableofcontents

\chapter{About SOSTOOLS v4.00}
The release of SOSTOOLS  v4.00 comes as we approach the twentieth anniversary of the original release of SOSTOOLS v1.00 back in April, 2002. SOSTOOLS was originally envisioned as a flexible tool for parsing and solving polynomial optimization problems, using the SOS tightening of polynomial positivity constraints, and capable of adapting to the ever-evolving fauna of applications of SOS. There are currently a variety of SOS programming parsers beyond SOSTOOLS, including YALMIP, Gloptipoly, SumOfSquares, and others. We hope SOSTOOLS remains the most intuitive, robust and adaptable toolbox for SOS programming. Recent progress in Semidefinite programming has opened up new possibilities for solving large Sum of Squares programming problems, and we hope the next decade will be one where SOS methods will find wide application in different areas. 

As many users of SOSTOOLS have already pointed out, parsing poses a significant computational and memory overhead in the process of solving a SOS problem, especially if that problem is large. Failure to solve large-scale SOS programming problems is many times traced back, not to the size of the resulting SDP, but rather to the computational and memory requirements of defining the SDP in the first place -- i.e. the parsing. In this sense, SOS parsers are not as efficient as they could be for  parsing large-scale SOS problems. In previous versions of SOSTOOLS we have already incorporated the multipoly toolbox -- stripping out much of the overhead associated with Matlab's symbolic math toolbox, and offering an alternative way of parsing SOS programs. The inclusion of the multipoly toolbox resulted in a significant decrease in parsing time and memory overhead. And yet, the computational and memory complexity of parsing large-scale SOS programming problems still remained a challenge, due to the tens of thousands of decision variables which appear in such large-scale SOS programming problems. 

In SOSTOOLS v4.00, we implement a parsing approach that reduces the computational and memory requirements of the parser below that of the SDP solver itself. To do this, we have completely re-developed the internal structure of our polynomial decision variables. Specifically, polynomial and SOS variable declarations made using \texttt{sossosvar}, \texttt{sospolyvar}, \texttt{sosmatrixvar}, etc now return a new polynomial structure, \texttt{dpvar}. This new polynomial structure, is optimized for decision variables, compatible with the multipoly toolbox, and documented in the enclosed \texttt{dpvar} guide, and isolates the scalar SDP decision variables in the SOS program from the independent variables used to construct the SOS program. As a result, the complexity of the parser scales almost linearly in the number of decision variables (significantly lower than the complexity of SDP solvers). While it is likely that most users will notice no difference in \textit{function} between SOSTOOLS 3.01 and SOSTOOLS 4.00, the internal parsing of these SOS programs has now been improved. As a result of these changes, almost all users will notice a significant increase in speed, with large-scale problems experiencing the most dramatic speedups. Parsing time is now always less than 10\% of time spent in the SDP solver and as a result, it is now possible to solve, e.g. local nonlinear stability analysis with 10 variables using degree 4 polynomials.

As mentioned above, the SOSTOOLS interface has changed very little with v4.00. Furthermore, we retain support for the use of symbolic variables (although we highly recommend using pvar). As a result, this user guide has changed very little. The only significant addition is support for the MOSEK~\cite{mosek_2021} SDP solver, the sospsimplify routine~\cite{sospsimplify}, and the addition of \texttt{sosquadvar} (which can be used for creating customized decision variables which fall outside the scope of \texttt{sossosvar}, \texttt{sospolyvar}, \texttt{sosmatrixvar}, etc.). 

The highlights of the latest SOSTOOLS release are listed below:
\begin{itemize}
\item Compatibility with newer versions of MATLAB including R2021a.
\item A new \texttt{dpvar} internal decision variable structure.
\item Interface to MOSEK~\cite{mosek_2021}.
\item Ability to call sospsimplify~\cite{sospsimplify}.
\item \texttt{sosquadvar} - A general-purpose function for declaring customized decision variables.
\end{itemize}

\chapter{Getting Started with \name}

\name\ is a free, third-party
MATLAB\footnote{A registered trademark of The MathWorks, Inc.}
toolbox for solving sum of squares programs. The techniques behind it
are based on the sum of squares decomposition for multivariate
polynomials \cite{ChoLR95}, which can be efficiently computed using
semidefinite programming \cite{VanB96}.  \name\ is developed as a
consequence of the recent interest in sum of squares polynomials
\cite{Par00,Par03,Sho87,ChoLR95,Rez00,Nesterovpoly,Lasserre}, partly due to the
fact that these techniques provide convex relaxations for many hard
problems such as global, constrained, and boolean optimization.

Besides the optimization problems mentioned above, sum of squares
polynomials (and hence \name) find applications in many other
areas. This includes control theory problems, such as:
search for Lyapunov functions to prove stability of a
dynamical system, computation of tight upper bounds for the structured
singular value $\mu$ \cite{Par00}, and stabilization of nonlinear
systems \cite{RanP00}. Some examples related to these problems, as
well as several other optimization-related examples, are provided and
solved in the demo files that are distributed with \name.

In the next two sections, we will provide a quick overview on sum of
squares polynomials and programs, and show the role of \name\ in sum
of squares programming.

\section{Sum of Squares Polynomials and Sum of Squares Programs}
A multivariate polynomial
$p(x_1,...,x_n)\triangleq p(x)$ is a sum of squares (SOS, for brevity), if there exist
polynomials $f_1(x),...,f_m(x)$ such that
\begin{equation}
p(x) = \sum_{i=1}^{m}f_i^2(x).
\label{eq:sospoly}
\end{equation}
It is clear that $f(x)$ being an SOS naturally implies $f(x)\geq 0$
for all $x \in \R^n$. For a (partial) converse statement, we remind
you of the equivalence, proven by Hilbert, between ``nonnegativity''
and ``sum of squares'' in the following cases:
\begin{itemize}
\item Univariate polynomials, any (even) degree.
\item Quadratic polynomials, in any number of variables.
\item Quartic polynomials in two variables.
\end{itemize}
(see \cite{Rez00} and the references therein). In the general
multivariate case, however, $f(x)\geq 0$ in the usual sense does not
necessarily imply that $f(x)$ is SOS. Notwithstanding this fact, the
crucial thing to keep in mind is that, while being stricter, the
condition that $f(x)$ is SOS is much more \emph{computationally
  tractable} than nonnegativity \cite{Par00}. At the same time,
practical experience indicates that replacing nonnegativity with the
SOS property in many cases leads to the exact solution.

The SOS condition (\ref{eq:sospoly}) is equivalent to the existence of
a positive semidefinite matrix $Q$, such that
\begin{equation}
p(x) = Z^T(x)QZ(x),
\label{Grameq}
\end{equation}
where $Z(x)$ is some properly chosen vector of monomials. Expressing an
SOS polynomial using a quadratic form as in (\ref{Grameq}) has also
been referred to as the Gram matrix method \cite{ChoLR95,PowW98}.

As hinted above, sums of squares techniques can be used to provide
tractable relaxations for many hard optimization problems.  A very
general and powerful relaxation methodology, introduced in
\cite{Par00,Par03}, is based on the \emph{Positivstellensatz}, a
central result in real algebraic geometry. Most examples in this
manual can be interpreted as special cases of the practical
application of this general relaxation method.
In this type of relaxations, we are interested in finding polynomials
$p_i(x)$, $i=1,2,...,\hat N$ and sums of squares $p_i(x)$ for $i=(\hat
N+1),...,N$ such that
\begin{eqnarray*}
& & a_{0,j}(x) + \sum_{i=1}^N p_i(x)a_{i,j}(x)
= 0,\text{\hspace{1em} for }j = 1,2,\ldots,J,
\end{eqnarray*}
where the $a_{i,j}(x)$'s are some given constant coefficient
polynomials. Problems of this type will be termed ``sum of squares
programs'' (SOSP).  Solutions to SOSPs like the above provide
certificates, or \emph{Positivstellensatz refutations}, which can be
used to prove the nonexistence of real solutions of systems of
polynomial equalities and inequalities (see \cite{Par03} for details).


The basic feasibility problem in SOS programming will be formulated as
follows: \myboxed{
\noindent \textbf{FEASIBILITY}:
\vspace{3 mm}\par
\noindent Find \vspace{2mm}\\
\hspace*{2em} polynomials $p_i(x)$, \hspace{2em} for $i=1,2,\ldots,\hat N$\\
\hspace*{2em} sums of squares $p_i(x)$, \hspace{2em} for $i=(\hat N+1),\ldots,N$\vspace{2mm}\\
such that
\begin{eqnarray}
& & a_{0,j}(x) + \sum_{i=1}^N p_i(x)a_{i,j}(x)
= 0,\text{\hspace{1em} for }j = 1,2,\ldots,\hat J, \label{Constraint1}\\
& & a_{0,j}(x) + \sum_{i=1}^N p_i(x)a_{i,j}(x)
\hspace{1em}\text{are sums of squares $(\geq 0)$\footnotemark,} \nonumber\\
& & \text{ \hspace{2em} for }j = (\hat J +
1),(\hat J + 2),\ldots, J. \label{Constraint2}
\end{eqnarray}}
\footnotetext{ Whenever constraint $f(x)\geq 0$ is encountered in
an SOSP, it should always be interpreted as ``$f(x)$ is an
SOS''.}

\noindent In this formulation, the $a_{i,j}(x)$ are
given scalar constant coefficient polynomials. The $p_i(x)$'s will be
termed \emph{SOSP variables}, and the constraints
(\ref{Constraint1})--(\ref{Constraint2}) are termed \emph{SOSP
  constraints}. The feasible set of this problem is convex, and as a
consequence SOS programming can in principle be solved using the
powerful tools of \emph{convex optimization}~\cite{BoyV04}.

It is obvious that the same program can be formulated in terms of
constraints (\ref{Constraint1}) only, by introducing some extra sums
of squares as slack program variables. However, we will keep this more
explicit notation for its added flexibility, since in most cases it
will help make the problem statement clearer.


Since many problems are more naturally formulated using inequalities,
we will call the constraints (\ref{Constraint2}) ``inequality
constraints'', and denote them by $\geq 0$. It is important, however,
to keep in mind the (possible) gap between nonnegativity and SOS.

Besides pure feasibility, the other natural class of problems in
convex SOS programming involves optimization of an objective function
that is linear in the coefficients of $p_i(x)$'s. The general form of
such optimization problem is as follows: \myboxed{
\noindent \textbf{OPTIMIZATION}:
\vspace{3 mm}\par
\noindent Minimize the linear objective function
\[
w^Tc,
\]
where $c$ is a vector formed from the (unknown) coefficients of
\vspace{2mm}\\
\hspace*{2em} polynomials $p_i(x)$, \hspace{2em} for $i=1,2,...,\hat N$\\
\hspace*{2em} sums of squares $p_i(x)$, \hspace{2em} for $i=(\hat N+1),...,N$\vspace{2mm}\\
such that
\begin{eqnarray}
& & a_{0,j}(x) + \sum_{i=1}^N p_i(x)a_{i,j}(x)
= 0,\text{\hspace{1em} for }j = 1,2,...,\hat J, \label{Constraint3}\\
& & a_{0,j}(x) + \sum_{i=1}^N p_i(x)a_{i,j}(x)
\hspace{1em}\text{are sums of squares $(\geq 0)$,} \nonumber\\
& & \text{ \hspace{2em} for }j = (\hat J +
1),(\hat J + 2),..., J, \label{Constraint4}
\end{eqnarray}}
\noindent In this formulation, $w$ is the vector of weighting coefficients
for the linear objective function.

Both the feasibility and optimization problems as formulated above are
quite general, and in specific cases reduce to well-known problems.
In particular, notice that if all the unknown polynomials $p_i$ are
restricted to be constants, and the $a_{i,j}, b_{i,j}$ are quadratic
forms, then we exactly recover the standard linear matrix inequality
(LMI) problem formulation.  The extra degrees of freedom in SOS
programming are actually a bit illusory, as every SOSP can be exactly
converted to an equivalent semidefinite program \cite{Par00}.
Nevertheless, for several reasons, the problem specification outlined
above has definite practical and methodological advantages, and
establishes a useful framework within which many specific problems can
be solved, as we will see later in Chapter~\ref{Applications}.


\section{What \name\ Does}

\begin{figure}[!tbp]
  \centering
  \includegraphics[height=.3\textheight]{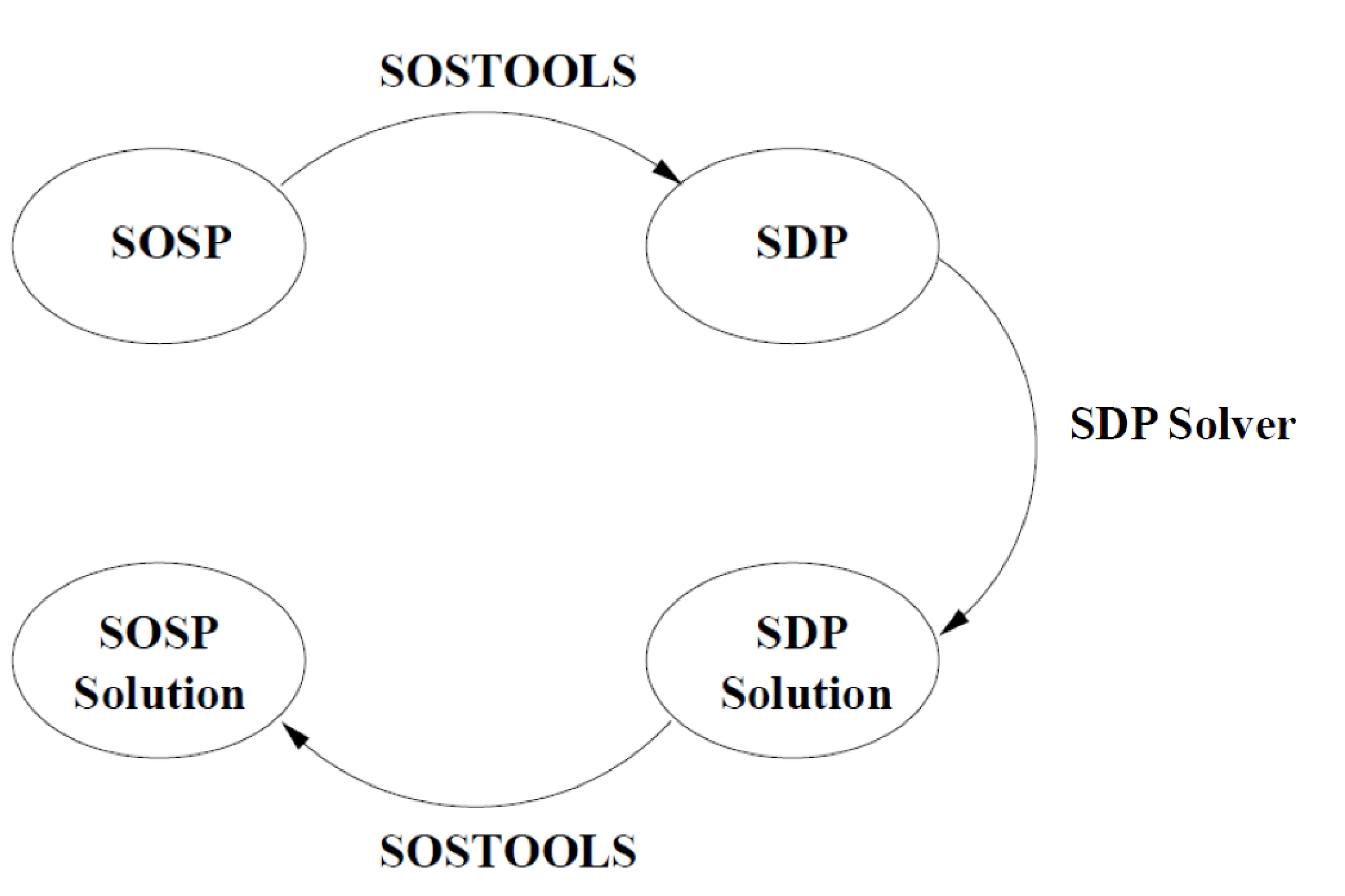}
  \caption{Diagram depicting relations between sum of squares program (SOSP), semidefinite program (SDP),
  \name, and SeDuMi/SDPT3/CSDP/SDPNAL/SDPNAL+/SDPA/Mosek.}\label{Fig1}
\end{figure}

Currently, sum of squares programs are solved by reformulating them as
semidefinite programs (SDPs), which in turn are solved efficiently
e.g.\ using interior point methods.  Several commercial as well as
non-commercial software packages are available for solving SDPs. While
the conversion from SOSPs to SDPs can be manually performed for small
size instances or tailored for specific problem classes, such a
conversion can be quite cumbersome to perform in general.  It is
therefore desirable to have a computational aid that automatically
performs this conversion for general SOSPs. This is exactly where
\name\ comes to play.  It automates the conversion from SOSP to SDP,
calls the SDP solver, and converts the SDP solution back to the
solution of the original SOSP.  At present there is an interface between \name\ and the following free MATLAB
based SDP solvers: i) SeDuMi \cite{Stu99},  ii) SDPT3 \cite{TohTT99}, iii) CSDP \cite{Bor99}, iv) SDPNAL \cite{ZhaST10}, v) SDPNAL+ \cite{YanST15}, vi) SDPA \cite{YFNNFKG10}, and MOSEK~\cite{mosek_2021}. This whole
process is depicted in Figure~\ref{Fig1}.

In the original release of \name, polynomials are implemented solely as symbolic objects, making
full use of the capabilities of the MATLAB Symbolic Math Toolbox. This
gives to the user the benefit of being able to do all polynomial
manipulations using the usual arithmetic operators: \verb"+",
\verb"-", \verb"*", \verb"/", \verb"^"; as well as operations such as
differentiation, integration,
point evaluation, etc. In addition, this provides the possibility of
interfacing with the Maple\footnote{A registered trademark of Waterloo
Maple Inc.} symbolic engine and the Maple library (which is very
advantageous). Unfortunately, not all users have access to the MATLAB Symbolic Math Toolbox. Furthermore, the Symbolic Math Toolbox has rather high computation and memory complexity. As a result, starting with SOSTOOLS v3.01 (although we continue to support the Symbolic Math Toolbox) we now highly recommend using our alternative custom-built polynomial multipoly (\texttt{pvar}) object to construct SOSP problems. Significantly, starting with SOSTOOLS v4.00, use of these multipoly objects allows the user to leverage the significant speedups associated with the internal \texttt{dpvar} decision variable structure.

The SOSTOOLS user interface has been designed to be as simple, as easy to use,
and as transparent as possible, while keeping a large degree of
flexibility. An SOSP is created by declaring SOSP variables (e.g., the
$p_i(x)$'s in Section~1.1), adding SOSP constraints, setting the
objective function, and so forth.  After the program is created,
one function is called to run the solver and finally the solutions to
SOSP are retrieved using another function.  These steps will be presented in
more details in Chapter~2.

Alternatively, ``customized'' functions for special problem classes
(such as Lyapunov function computation, etc.) can be directly used,
with no user programming whatsoever required. These are presented in
the first three sections of Chapter~3.

\section{System Requirements and Installation Instruction}
To install and run SOSTOOLS v4.00, MATLAB R2009a or later is required. Older versions of both MATLAB and the symbolic toolbox (if using \texttt{sym} polynomial objects) may be sufficient however SOSTOOLS v4.00 has only been tested on versions 2009a -- 2021a. Here is a list of requirements:
\begin{itemize}
\item MATLAB R2009a or later.

\item Symbolic Math Toolbox version 5.7 or later (if using \texttt{sym} polynomial objects).
\item A PC with minimum 8GB RAM (more is better)

\item One of the following SDP solvers: SeDuMi, SDPT3, CSDP,  SDPNAL, SDPA, and MOSEK. Each solver must be installed before SOSTOOLS can be used. The user is referred to the relevant documentation to see how this is done\footnote{It is not recommended that both SDPNAL and SDPT3 be on the MATLAB path simultaneously.}\footnote{To use the CSDP solver please note that the CSDP binary file must be in the working directory and not simply be on the MATLAB path.}. The solvers can be downloaded from:

SeDuMi: \verb"http://sedumi.ie.lehigh.edu"

SDPT3: \verb"http://www.math.nus.edu.sg/~mattohkc/sdpt3.html"

CSDP: \verb"https://projects.coin-or.org/Csdp/"

SDPNAL: \verb"http://www.math.nus.edu.sg/~mattohkc/SDPNAL.html"

SDPNAL+: \verb"http://www.math.nus.edu.sg/~mattohkc/SDPNALplus.html"

SDPA: \verb"http://sdpa.sourceforge.net/index.html"

MOSEK: \verb"https://www.mosek.com/downloads/"

\end{itemize}

Note that if you do not have access to the Symbolic Toolbox then SOSTOOLS v4.00 can be used with the multivariate polynomial toolbox and any version of MATLAB.

\name\ can be easily run on a
UNIX workstation, on  Windows operating systems, and Mac OSX. It
utilizes MATLAB sparse matrix representation for good performance and
to reduce the amount of memory needed. To give an illustrative figure
of the computational load, all examples in Chapter~\ref{Applications} except
the $\mu$ upper bound example, are solved in less than 10 seconds by
\name\ running on a PC with Intel Celeron 700 MHz processor and 96
MBytes of RAM. Even the $\mu$ upper bound example is solved in less
than 25 seconds using the same system.

\name\ is available for free under the GNU General Public License. The most recent version of SOSTOOLS can be downloaded from GitHub at \verb"https://github.com/oxfordcontrol/SOSTOOLS". Previous versions can be downloaded from
\verb"http://www.eng.ox.ac.uk/control/sostools/" \linebreak or
\verb"http://www.cds.caltech.edu/sostools" or
\verb"http://www.mit.edu/~parrilo/sostools/" or \verb"http://control.asu.ed/sostools/".  Once you download
the zip file, you should extract its contents to the directory where
you want to install \name. In UNIX, you may use
\begin{matlab}
unzip -U \name.nnn.zip -d \verb"your_dir"
\end{matlab}
where \verb"nnn" is the version number, and \verb"your_dir" should be replaced
by the directory of your choice. In Windows operating systems, you may use
programs like Winzip to extract the files.

After this has been done, you must add the \name\ directory and its subdirectories to the MATLAB path.
This is done in MATLAB by choosing the menus File \verb"-->" Set Path \verb"-->" Add with Subfolders ..., and
then typing the name of \name\  main directory. This completes the \name\ installation. Alternatively, run the script \verb"addsostools.m" from the SOSTOOLS  directory.

\section{Other Things You Need to Know}
The directory in which you install \name\ contains several subdirectories.
Two of them are:
\begin{itemize}
\item \verb"sostools/docs" : containing this user's manual and license file.
\item \verb"sostools/demos" : containing several demo files.
\end{itemize}
The demo files in the second subdirectory above implement the SOSPs corresponding to
examples in Chapter~\ref{Applications}.

Throughout this user's manual, we use the \verb"typewriter" typeface to
denote MATLAB variables and functions, MATLAB commands that you should type, and results given by MATLAB.
MATLAB commands that you should type will also be denoted by the symbol \verb">>" before the commands.
For example,
\begin{matlab}
\begin{verbatim}
>> x = sin(1)

x =

    0.8415

\end{verbatim}
\end{matlab}
In this case, \verb"x = sin(1)" is the command that you type, and \verb"x = 0.8415"
is the result given by MATLAB.

Finally, you can send bug reports, comments, and suggestions to
\verb"sostools@cds.caltech.edu", or directly on GitHub. Any feedback is greatly appreciated.




\chapter{Solving Sum of Squares Programs}\label{ch:SOSP}

\name\ can solve two kinds of sum of squares programs: the feasibility and
optimization problems, as formulated in Chapter~2.  To define and
solve an SOSP using \name, you simply need to follow
these steps:
\begin{enumerate}
\item Initialize the SOSP.
\item Declare the SOSP variables.
\item Define the SOSP constraints.
\item Set objective function (for optimization problems).
\item Call solver.
\item Get solutions.
\end{enumerate}
In the next sections, we will describe each of these steps in detail. But first,
we will look at how polynomials are represented and manipulated in \name.

\section{Polynomial Representation and Manipulations}
\label{Polynomial Representation}

Polynomials in \name\ can have representation as symbolic objects, using the MATLAB Symbolic Toolbox.
Typically, a polynomial is created by first declaring its independent variables and then constructing
it using the usual algebraic manipulations. For example, to create a polynomial
$p(x,y) = 2x^2+3xy+4y^4$,
you declare the independent variables $x$ and $y$ by typing
\begin{matlab}
\begin{verbatim}
>> syms x y;
\end{verbatim}
\end{matlab}
and then construct $p(x,y)$ as follows:
\begin{matlab}
\begin{verbatim}
>> p = 2*x^2 + 3*x*y + 4*y^4

p =

2*x^2+3*x*y+4*y^4
\end{verbatim}
\end{matlab}

Polynomials such as the one created above can then be manipulated using the usual operators:
\verb"+",
\verb"-", \verb"*", \verb"/", \verb"^". Another operation which is particularly useful for control-related
problems such as Lyapunov function search is differentiation, which can be done using the function \verb"diff".
For instance, to find the partial derivative $\frac{\partial p}{\partial x}$, you should type
\begin{matlab}
\begin{verbatim}
>> dpdx = diff(p,x)

dpdx =

4*x+3*y

\end{verbatim}
\end{matlab}
For other types of symbolic manipulations, we refer you to the manual and help comments of the Symbolic Math Toolbox.

Starting with SOSTOOLS v3.01, users are encouraged to use the alternative custom-built multipoly polynomial object in the Multivariate Polynomial Toolbox -- a freely available toolbox for constructing and manipulating
multivariate polynomials included in SOSTOOLS directory \texttt{multipoly}. In the remainder of the section, we give a
brief introduction to the multipoly polynomial objects in SOSTOOLS.

Independent polynomial variables are created with the \texttt{pvar} command.  For example, the following command creates three indpendent variables:
\begin{quote}
\begin{verbatim}
>> pvar x1 x2 x3
\end{verbatim}
\end{quote}
New polynomial objects can now be created from these variables, and
manipulated using, e.g. standard addition, multiplication, and integer
exponentiation functions:
\begin{quote}
\begin{verbatim}
>> p = x3^4+5*x2+x1^2
p =
  x3^4 + 5*x2 + x1^2
\end{verbatim}
\end{quote}
Matrices of polynomials can be created from polynomials using
horizontal/vertical concatenation and block diagonal augmentation, e.g.:
\begin{quote}
\begin{verbatim}
>> M1 = blkdiag(p,2*x2)
M1 =
  [ x3^4 + 5*x2 + x1^2 ,    0 ]
  [                  0 , 2*x2 ]
\end{verbatim}
\end{quote}
Naturally, it is also possible to build new polynomial matrices from
already constructed submatrices.
Elements of a polynomial matrix can be referenced and assigned using
the standard MATLAB referencing notation:
\begin{quote}
\begin{verbatim}
>> M1(1,2)=x1*x2
M1 =
  [ x3^4 + 5*x2 + x1^2 , x1*x2 ]
  [                  0 ,  2*x2 ]
\end{verbatim}
\end{quote}

The internal data structure for an $N\times M$ polynomial matrix of
$V$ variables and $T$ terms consists of a $T \times NM$ sparse
coefficient matrix, a $T \times V$ degree matrix, and a $V \times 1$
cell array of variable names.  This information can be easily accessed
through the MATLAB field accessing operators: \texttt{p.coefficient},
\texttt{p.degmat}, and \texttt{p.varname}. The access to fields uses a
case insensitive, partial-match.  Thus abbreviations, such as
\texttt{p.coef}, can also be used to obtain the coefficients, degrees,
and variable names.  Many additional operations exist in the multipoly toolbox such as trace, transpose, determinant,
differentiation, logical equal/not equal etc. In addition, multipoly includes converters between the symbolic and multipoly formats using \texttt{s2p} and \texttt{p2s}. Finding the roots of a multipoly polynomial can we performed using \texttt{psolve}. See the included multipoly documentation for a complete list of functions.

The input to the SOSTOOLS commands can be specified using either the
symbolic objects or multipoly polynomial objects
(although they cannot be mixed). The advantage of using the multipoly polynomial format is that it unlocks the internal use of the \texttt{dpvar} structure for polynomial decision variables as defined in Section~\ref{sec:decvars}. 

The \texttt{dpvar} polynomial format is similar to the multipoly polynomial format with the exception that it includes a list of decision variables to be used in the SDP solver. \texttt{dpvar} polynomials are only created using decision variable declarations as described in Section~\ref{sec:decvars}. Once a \texttt{dpvar} polynomial has been created, it can be manipulated in the same manner as a symbolic or multipoly object, with the exception that two \texttt{dpvar} objects cannot be multiplied, as this would create a bilinearity in the resulting SDP. \texttt{dpvar} structures can, however, be combined with multipoly objects using addition, subtraction, and multiplication.

The internal data structure for an $M\times N$ polynomial matrix of
$V$ independent variables, $D$ decision variables and $T$ monomials consists of a $(D+1)M \times NT$ sparse coefficient matrix, a $T \times V$ degree matrix, a $V \times 1$ cell array of variable names, and a $D \times 1$ cell array of decision variable names.  This information can be easily accessed
through the MATLAB field accessing operators: \texttt{p.coefficient},
\texttt{p.degmat}, \texttt{p.varname}, and \texttt{p.dvarname}. See the dpvar user manual included in the SOSTOOLS documentation for additional details.


\section{Initializing a Sum of Squares Program}
A sum of squares program is initialized using the command \verb"sosprogram". A vector containing independent
variables in the program has to be given as an argument to this function. Thus, if the polynomials
in our program have $x$ and $y$ as the independent variables, then we initialize the SOSP using
\begin{matlab}
\begin{verbatim}
>> Program1 = sosprogram([x;y]);
\end{verbatim}
\end{matlab}
The command above will initialize an empty SOSP called \verb"Program1".

\paragraph{Symbolic or Multipoly Format} When \texttt{sosprogam} is called and \texttt{Program1} initialized, SOSTOOLS will detect whether the variables $x$ and $y$ are symbolic or multipoly format. Once \texttt{sosprogram} has been called, this polynomial type is hardcoded into \texttt{Program1} and cannot be changed. 

\section{Variable Declaration}\label{sec:decvars}
After the program is initialized, the SOSP decision variables have to
be declared. There are five functions used for this purpose,
corresponding to variables of these types:
\begin{itemize}
\item Scalar decision variables.
\item Polynomial variables.
\item Sum of squares variables.
\item Matrix of polynomial variables.
\item Matrix of sum of squares variables.
\item Customized variables.
\end{itemize}
Each of them will be described in the following subsections.

\subsection{Scalar Decision Variables}
Scalar decision variables in an SOSP are meant to be unknown
scalar constants. The variable $\gamma$ in \verb"sosdemo3.m" (see
Section 3.3) is an example of such a variable. These variables can be
declared either by specifying them when an SOSP is initialized with
\verb"sosprogram", or by declaring them later using the function
\verb"sosdecvar".

\paragraph{Symbolic Math Format} To declare decision variables using the Symbolic Math Toolbox, you must first create symbolic objects corresponding to your decision variables.  This is performed using the functions \verb"syms" or \verb"sym" from the Symbolic Math Toolbox, in a way similar to the one you use to define independent variables in Section~\ref{Polynomial Representation}.  As explained earlier, you can declare the decision variables when you initialize an SOSP, by
giving them as a second argument to \verb"sosprogram". Thus, to
declare variables named \verb"a" and \verb"b", use the following
command:
\begin{matlab}
\begin{verbatim}
>> syms x y a b;
>> Program2 = sosprogram([x;y],[a;b]);
\end{verbatim}
\end{matlab}
Alternatively, you may declare these variables after the SOSP is
initialized, or add some other decision variables to the program,
using the function \verb"sosdecvar". For example, the sequence of
commands above is equivalent to
\begin{matlab}
\begin{verbatim}
>> syms x y a b;
>> Program3 = sosprogram([x;y]);
>> Program3 = sosdecvar(Program3,[a;b]);
\end{verbatim}
\end{matlab}
and also equivalent to
\begin{matlab}
\begin{verbatim}
>> syms x y a b;
>> Program4 = sosprogram([x;y],a);
>> Program4 = sosdecvar(Program4,b);
\end{verbatim}
\end{matlab}
\paragraph{Multipoly Format}
When using multipoly objects, the process is slightly different that for symbolic variables. Specifically, decision variables should be created using the command \texttt{dpvar} and independent variables should be created using \texttt{pvar}. Considering the syntax used above, we would have either
\begin{matlab}
\begin{verbatim}
>> pvar x y;
>> dpvar a b;
>> Program2 = sosprogram([x;y],[a;b]);
\end{verbatim}
\end{matlab}
or 
\begin{matlab}
\begin{verbatim}
>> pvar x y;
>> dpvar a b;
>> Program3 = sosprogram([x;y]);
>> Program3 = sosdecvar(Program3,[a;b]);
\end{verbatim}
\end{matlab}
Note that for compatibility with legacy codes, the user may still use 
\begin{matlab}
\begin{verbatim}
>> pvar x y a b;
>> Program2 = sosprogram([x;y],[a;b]);
\end{verbatim}
\end{matlab}
This legacy usage will result in a \textit{very} slight reduction in efficiency of the parser.

\subsection{Scalar Polynomial Variables}

Polynomial variables in a sum of squares program are simply polynomials with unknown coefficients (e.g.\ $p_1(x)$ in
the feasibility problem formulated in Chapter~1). Polynomial
variables can obviously be created by declaring its unknown coefficients as decision variables, and then
constructing the polynomial itself via some algebraic manipulations. For example, to create a polynomial
variable $v(x,y) = ax^2+bxy+cy^2$, where $a$, $b$, and $c$ are the unknowns, you can use the following commands:
\begin{matlab}
\begin{verbatim}
>> Program5 = sosdecvar(Program5,[a;b;c]);
>> v = a*x^2 + b*x*y + c*y^2;
\end{verbatim}
\end{matlab}
However, such an approach would be inefficient for polynomials with
many coefficients.  In such a case, you should use the function
\verb"sospolyvar" to declare a polynomial variable:
\begin{matlab}
\begin{verbatim}
>> [Program6,v] = sospolyvar(Program6,[x^2; x*y; y^2]);
\end{verbatim}
\end{matlab}
In this case \verb"v" will be
\begin{matlab}
\begin{verbatim}
>> v

v =

coeff_1*x^2+coeff_2*x*y+coeff_3*y^2

\end{verbatim}
\end{matlab}
We see that \verb"sospolyvar" automatically creates
decision variables corresponding to monomials in the vector which is given as the second input argument to
it, and then constructs a polynomial variable from these coefficients and monomials.
This polynomial variable is returned as the second output argument of \verb"sospolyvar".
\vspace{2mm}

\noindent \textbf{NOTE:}
\begin{enumerate}
\item \verb"sospolyvar" and \verb"sossosvar" (see Section 2.3.4) name the unknown coefficients (scalar decision variables) in
a polynomial/SOS variable
by \verb"coeff_nnn", where \verb"nnn" is a number.
Names that begin with \verb"coeff_" are reserved for this purpose, and therefore must not be used elsewhere.

\item
By default, the decision variables \verb"coeff_nnn" created by \verb"sospolyvar" or \verb"sossosvar"
will only be available in the
function workspace, and therefore cannot be manipulated in the MATLAB workspace. Sometimes it is desirable to have
these decision variables available in the MATLAB workspace, such as
when we want to set an objective function of an SOSP that involves one or more of these variables.
In this case, a third argument \verb"'wscoeff'" has to be given to \verb"sospolyvar" or \verb"sossosvar".
For example, using
\begin{matlab}
\begin{verbatim}
>> [Program7,v] = sospolyvar(Program7,[x^2; x*y; y^2],'wscoeff');
>> v

v =

coeff_1*x^2+coeff_2*x*y+coeff_3*y^2

\end{verbatim}
\end{matlab}
you will be able to directly use \verb"coeff_1" and \verb"coeff_2" in the MATLAB workspace, as shown below.
\begin{matlab}
\begin{verbatim}
>> w = coeff_1+coeff_2

w =

coeff_1+coeff_2
\end{verbatim}
\end{matlab}\vspace{-1em}

\item \name\ requires monomials that are given as the second input argument to \verb"sospolyvar" and \verb"sossosvar"
to be unique, meaning that there are no repeated monomials.
\end{enumerate}

\subsection{An Aside: Constructing Vectors of Monomials}

We have seen in the previous subsection that for declaring SOSP variables using \verb"sospolyvar" we need to construct a
vector whose entries are monomials. While this can be done by creating the individual monomials and arranging
them as a vector, \name\ also provides a function, named \verb"monomials",
that can be used to construct a column vector of monomials
with some pre-specified degrees. This will be particularly useful when the vector contains a lot of monomials.
The function takes two arguments: the first argument
is a vector containing all independent variables in the monomials, and the second argument is
a vector whose entries are the degrees of monomials
that you want to create. As an example, to construct a vector containing all monomials in $x$ and $y$ of degree 1, 2,
and 3, type the following command:

\begin{matlab}
\begin{verbatim}
>> VEC = monomials([x; y],[1 2 3])
\end{verbatim}
\end{matlab}

\begin{matlab}
\begin{verbatim}

VEC =

[     x]
[     y]
[   x^2]
[   x*y]
[   y^2]
[   x^3]
[ x^2*y]
[ x*y^2]
[   y^3]

\end{verbatim}
\end{matlab}
We clearly see that \verb"VEC" contains all monomials in $x$ and $y$ of degree 1, 2, and 3.

For some problems, such as Lyapunov
stability analysis for linear systems with parametric uncertainty, it is
desirable to declare polynomials with a certain structure called the
\emph{multipartite} structure. See Section~\ref{Sec:Multipartite} for a more
thorough discussion on this kind of structure. Multipartite polynomials
are declared using a monomial vector that also has multipartite structure.
To construct multipartite monomial vectors, the command \verb"mpmonomials" can be used.
For example,
\begin{matlab}
\begin{verbatim}
>> VEC = mpmonomials({[x1; x2],[y1; y2],[z1]},{1:2,1,3})

VEC =

[    z1^3*x1*y1]
[    z1^3*x2*y1]
[  z1^3*x1^2*y1]
[ z1^3*x1*x2*y1]
[  z1^3*x2^2*y1]
[    z1^3*x1*y2]
[    z1^3*x2*y2]
[  z1^3*x1^2*y2]
[ z1^3*x1*x2*y2]
[  z1^3*x2^2*y2]

\end{verbatim}
\end{matlab}
will create a vector of multipartite monomials where the partitions of the independent variables
are $S_1=\{x_1,x_2\}$, $S_2=\{y_1,y_2\}$, and $S_3=\{z_1\}$, whose corresponding degrees
are 1--2, 1, and~3.

\subsection{Sum of Squares Variables}
 
Sum of squares variables are also polynomials with unknown coefficients, similar to polynomial variables
described in Section 3.3.2. The difference is, as its name suggests, that an SOS variable is constrained to be
an SOS. This
is imposed by internally representing an SOS variable in the Gram matrix form (cf.~Section 2.1),
\begin{equation}
p(x) = Z^T(x)QZ(x)\label{sossosvar}
\end{equation}
and requiring the coefficient matrix Q to be positive semidefinite.

To declare an SOS variable, you must use the function \verb"sossosvar". The monomial vector $Z(x)$ in
(\ref{sossosvar}) has to be given as the second input argument to the function. Like \verb"sospolyvar", this
function will automatically declare all decision variables corresponding to the matrix \verb"Q". For example,
to declare an SOS variable
\begin{equation}
p(x,y) = \begin{bmatrix}x\\ y\end{bmatrix}^T Q \begin{bmatrix}x\\ y\end{bmatrix},
\end{equation}
type
\begin{matlab}
\begin{verbatim}
>> [Program8,p] = sossosvar(Program8,[x; y]);
\end{verbatim}
\end{matlab}
where the second output argument is the name of the variable.
In this example, the coefficient matrix
\begin{equation}
Q = \begin{bmatrix}\verb"coeff_1" & \verb"coeff_3"\\ \verb"coeff_2" & \verb"coeff_4" \end{bmatrix} \label{Qexample}
\end{equation}
will be created by the function. When this matrix is substituted into the expression for $p(x,y)$, we obtain
\begin{equation}
p(x,y) = \verb"coeff_1"x^2 + (\verb"coeff_2"+\verb"coeff_3")xy + \verb"coeff_4"y^2,
\end{equation}
which is exactly what \verb"sossosvar" returns:
\begin{matlab}
\begin{verbatim}
>> p

p =

coeff_4*y^2+(coeff_2+coeff_3)*x*y+coeff_1*x^2
\end{verbatim}
\end{matlab}

We would like to note that at first the coefficient matrix does not appear to be symmetric, especially
because the number of decision variables (which seem to be independent) is the same as the number of
entries in the coefficient matrix.
However, some constraints are internally
imposed by the semidefinite programming solver SeDuMi/SDPT3 (which are used by \name) on some of these decision variables,
such that
the solution matrix obtained by the solver will be symmetric.
The primal formulation of a semidefinite program in
SeDuMi/SDPT3 uses $n^2$ decision variables to represent an $n\times n$ positive
semidefinite matrix, which is the reason why \name\ also uses $n^2$ decision variables for its $n \times n$
coefficient matrices.

SOSTOOLS includes a custom function \verb"findsos" that will compute, if feasible, the sum of squares decomposition of a polynomial $p(x)$ into the sum of $m$ polynomials $f_{i}^{2}(x)$ as in (\ref{eq:sospoly}), the Gram matrix \verb"Q" and vector of monomials  \verb"Z" corresponding to (\ref{sossosvar}). The function is called as shown below:
\begin{matlab}
\begin{verbatim}
>> [Q,Z,f] = findsos(P);
\end{verbatim}
\end{matlab}
where \verb"f" is a vector of length $m=rank(Q)$ containing the functions $f_{i}$. If the problem is infeasible then empty matrices are returned. This example is expanded upon in \verb"SOSDEMO1" in Chapter 4.


\subsection{Matrix Variables}

For many problems it may be necessary to construct  matrices of polynomial or sum of squares polynomials decision variables, i.e. matrices whose elements are themselves   polynomials or sum of squares polynomials with unknown coefficients. Such decision variables can be respectively declared using the \verb"sospolymatrixvar" or the \verb"sossosmatrixvar" function. The \verb"sospolymatrixvar" or the \verb"sossosmatrixvar" functions take three compulsory input arguments and an optional fourth symmetry argument. The first two arguments are of the same form as \verb"sospolyvar" and \verb"sossosvar", the first being the sum of squares program, \verb"prog", and the second the vector of monomials $Z(x)$. The third argument is a row vector specifying the dimension of the matrix. We now illustrate a few simple examples of the use of the  \verb"sospolymatrixvar" function. First a SOSP must be initialized:
\begin{matlab}
\begin{verbatim}
>> syms x1 x2;
>> x = [x1 x2].';
>> prog = sosprogram(x);
\end{verbatim}
\end{matlab}
We will now declare two matrices \verb"P" and \verb"Q" both of dimension $2\times 2$ where the entries are real scalars, i.e. a degree 0 polynomial matrix. Furthermore we will add the constraint that \verb"Q" must be symmetric:
\begin{matlab}
\begin{verbatim}
>> [prog,P] = sospolymatrixvar(prog,monomials(x,0),[2 2]);
>> [prog,Q] = sospolymatrixvar(prog,monomials(x,0),[2 2],'symmetric');
>> P

    P =

    [coeff_1, coeff_2]
    [coeff_3, coeff_4]

>> Q

     Q =

    [coeff_5, coeff_6]
    [coeff_6, coeff_7]
\end{verbatim}
\end{matlab}
 To declare a symmetric matrix where the elements are homogenous quadratic polynomials the function \verb"sospolymatrixvar" is called with the following arguments:
\begin{matlab}
\begin{verbatim}
>> [prog,R] = sospolymatrixvar(prog,monomials(x,2),[2 2]);
>> R(1,1)

    ans =

    coeff_8*x1^2 + coeff_9*x1*x2 + coeff_10*x2^2

>> R(1,2)

   ans =

   coeff_11*x1^2 + coeff_12*x1*x2 + coeff_13*x2^2
\end{verbatim}
\end{matlab}
\begin{matlab}
\begin{verbatim}
>> R(2,1)

     ans

    coeff_11*x1^2 + coeff_12*x1*x2 + coeff_13*x2^2

>> R(2,2)

   ans =

    coeff_14*x1^2 + coeff_15*x1*x2 + coeff_16*x2^2
\end{verbatim}
\end{matlab}
In the next section it will be shown how these matrix variables can be incorporated as constraints in sum of squares optimization problems.

\subsection{Customized Variables} 

Occasionally, the user may wish to define a decision variable which is not included in the list of standard decision variable types. When using the multipoly polynomial format, this can be achieved using \texttt{sosquadvar}. 

\begin{matlab}
\begin{verbatim}
>> [Program9,p] = sosquadvar(Program9,Z1,Z2,m,n,option)
\end{verbatim}
\end{matlab}
creates a variable of the form
\begin{align*}
P&=\left(I_{m} \otimes Z_{1}(x)\right)^TQ\;\left(I_{n} \otimes Z_{2}(y)\right)\\ \small
&=\begin{bmatrix}
Z_{1}(x) &&\\
&\hspace{-8mm}\ddots&\\
&&\hspace{-8mm}Z_{1}(x)
\end{bmatrix}^TQ\begin{bmatrix}
Z_{2}(y)&&\\
&\hspace{-8mm}\ddots&\\
&&\hspace{-8mm}Z_{2}(y)
\end{bmatrix}
\end{align*}
The matrix $Q$ is the matrix of decision variables. The \texttt{option} specifies if the matrix $Q$ should be constrained to be positive semidefinite or symmetric or neither. This class of variables encompasses existing classes of decision variables and is especially useful for positive kernel matrices. The inputs to \texttt{sosquadvar} are
\begin{itemize}
\item Z1 - A column vector of monomials (pvar) to be multiplied on the left
\item Z2 - A column vector of monomials (pvar) to be multiplied on the right
\item m - a scalar indicating the row dimension of the output. If empty, it defaults to 1
\item n - a scalar indicating the column dimension of the output. If empty, it defaults to 1
\item option - 'pos' if Q should be a positive semidefinite matrix 'sym' if Q should be a symmetric matrix (note this may not result in P being symmetric). In order to use the positive or symmetric options, the length of the monomial vectors Z1 and Z2 should be identical. In addition, for this case, we require $m$=$n$.
\end{itemize}
To illustrate, 
\begin{matlab}
\begin{verbatim}
>> pvar x y z
>> [Program9,p] = sosquadvar(Program9,[x y],[x z],1,2,'pos')
\end{verbatim}
\end{matlab}
creates a variable of the form
\begin{align*}
P&=
\begin{bmatrix}
x\\
y
\end{bmatrix}^TQ
\begin{bmatrix}
x&0\\
z&0\\
0&x\\
0&z
\end{bmatrix}
\end{align*}
with $Q \ge 0$ a positive semidefinite matrix.

\paragraph{Cellular Inputs and Outputs} \texttt{sosquadvar} can also take cellular arguments in order to allow for variable structures where positivity is coupled between several polynomial decision variables. 
\begin{matlab}
\begin{verbatim}
>> [Program9,P]=sosquadvar(Program9,Z1c,Z2c,m,n,option)
\end{verbatim}
\end{matlab}
In this case, we allow for a multipartite structure. Z1c and Z2c are cells of monomial column vectors. The number of cells in Z1c and Z2c need not be the same unless the 'pos' or 'sym' options are used. Note that m and n are now expected to be vectors with length matching the number of cells in Z1c and Z2c, respectively and indicating the number of rows and columns to be output for each of the cells in Z1c and Z2c. The ouput in this case is a cellular structure, $P$, where $P\{i,j\} \in \R^{n_i \times n_j}$ has the form
\[
P\{i,j\}(x,y)=(I_{m_i} \otimes Z1c\{i\}(x))^TQ_{ij} (I_{n_j}\otimes Z2c\{j\}(y))
\]
where $Q_{ij}$ is an dimension-appropriate block partition of the decision variable $Q$ where the if the 'pos' option is selected, $Q\ge 0$. Note that this 'pos' option does \textbf{not} imply the sub-blocks $Q_{ij}$ are positive semidefinite for $i \neq j$.

\section{Adding Constraints}
Sum of squares program constraints such as (\ref{Constraint1})--(\ref{Constraint2}) are
added to a sum of squares program using the functions \verb"soseq" and \verb"sosineq".

\subsection{Equality Constraints}
For adding an equality constraint to a sum of squares program, you must use the function \verb"soseq".
As an example, assume that \verb"p" is an SOSP variable, then
\begin{matlab}
\begin{verbatim}
>> Program9 = soseq(Program9,diff(p,x)-x^2);
\end{verbatim}
\end{matlab}
will add the equality constraint
\begin{equation}
\frac{\partial p}{\partial x} - x^2 = 0
\end{equation}
to \verb"Program9".

\subsection{Inequality Constraints}
Inequality constraints are declared using the function \verb"sosineq", whose basic syntax is
similar to \verb"soseq". For example, type
\begin{matlab}
\begin{verbatim}
>> Program1 = sosineq(Program1,diff(p,x)-x^2);
\end{verbatim}
\end{matlab}
to add the inequality constraint\footnotemark
\begin{equation}
\frac{\partial p}{\partial x} - x^2 \geq 0.
\end{equation}
\footnotetext{We remind you that $\frac{\partial p}{\partial x} - x^2 \geq 0$ has to be interpreted as
$\frac{\partial p}{\partial x} - x^2$ being a sum of squares. See the discussion in Section 1.1.}
\addtocounter{footnote}{1}%

\noindent However, several differences do exist. In particular, a third argument can be given to \verb"sosineq" to
handle the following cases:
\begin{itemize}
\item When there is only one independent variable in the SOSP (i.e., if the polynomials
are univariate), a third argument can be given to specify the range of independent variable for which the inequality
constraint has to be satisfied. For instance, assume that \verb"p" and \verb"Program2" are respectively univariate
polynomial and univariate SOSP, then
\begin{matlab}
\begin{verbatim}
>> Program2 = sosineq(Program2,diff(p,x)-x^2,[-1 2]);
\end{verbatim}
\end{matlab}
will add the constraint
\begin{equation}
\frac{\partial p}{\partial x} - x^2 \geq 0,\hspace{1em}\text{for }-1\leq x \leq 2
\end{equation}
to the SOSP. See Sections 3.7 and 3.8 for application examples where this option is used.

\item When the left side of the inequality is a high degree sparse
  polynomial (i.e., containing a few nonzero terms), it is
  computationally more efficient to impose the SOS condition using a
  reduced set of monomials (see \cite{Rez78}) in the Gram matrix form.
  This will result in a smaller size semidefinite program, which is
  easier to solve.  By default, \name\ does not try to obtain this
  optimal reduced set of monomials, since this itself takes an
  additional amount of computational effort (however, \name\ always does
  some reasonably efficient and computationally cheap heuristics to reduce the set of monomials).
  The optimal reduced
  set of monomials will be computed and used only if a third argument
  \verb"'sparse'" is given to \verb"sosineq", as illustrated by the
  following command,
\begin{matlab}
\begin{verbatim}
>> Program3 = sosineq(Program3,x^16+2*x^8*y^2+y^4,'sparse');
\end{verbatim}
\end{matlab}
which tests whether or not $x^{16}+2x^8y^2+y^4$ is a sum of squares. See Section~\ref{Sec:Sparse} for a discussion
on exploiting sparsity.

\item A special sparsity structure that can be easily handled is the multipartite structure.
When a polynomial has this kind of structure, the optimal reduced set of monomials
in $Z^T(x)QZ(x)$ can be obtained with a low computational effort. For this, however, it is necessary
to give a third argument \verb"'sparsemultipartite'" to \verb"sosineq", as well as the partition of
the independent variables which form the multipartite structure. As an example,
\begin{matlab}
\begin{verbatim}
>> p = x1^4*y1^2+2*x1^2*x2^2*y1^2+x2^2*y1^2;
>> Program3 = sosineq(Program3,p,'sparsemultipartite',{[x1,x2],[y1]});
\end{verbatim}
\end{matlab}
tests whether or not the multipartite (corresponding to partitioning the independent variables to $\{x_1,x_2\}$ and $\{y_1\}$)
polynomial $x_1^4y_1^2+2x_1^2x_2^2y_1^2+x_2^2y_1^2$
is a sum of squares. See Section~\ref{Sec:Multipartite} for a discussion
on the multipartite structure.

\end{itemize}
\textbf{NOTE:} The function \verb"sosineq" will accept matrix arguments in addition to scalar arguments. Matrix arguments are \emph{not} treated element-wise inequalities. To avoid confusion  it is suggested that \verb"sosmatrixineq" is used when defining matrix inequalities.

\subsection{Exploiting Sparsity}
\label{Sec:Sparse}

For a polynomial $p(x)$, the complexity of computing the sum of
squares decomposition $p(x)=\sum_i p_i^2(x)$ (or equivalently,
$p(x)=Z(x)^TQZ(x)$, where $Z(x)$ is a vector of monomials ---
see~\cite{Par00} for details) depends on two factors: the number of
variables and the degree of the polynomial. However when $p(x)$ has
special structural properties, the computation effort can be notably
simplified through the reduction of the size of the semidefinite
program, removal of degeneracies, and better numerical
conditioning. Since the initial version of SOSTOOLS, Newton polytopes
techniques have been available via the optional argument
\texttt{'sparse'} to the function \texttt{sosineq}.

The notion of sparseness for multivariate polynomials is
stronger than the one commonly used for matrices. While in the matrix
case this word usually means that many coefficients are zero, in the
polynomial case the specific vanishing pattern is also taken into
account. This idea is formalized by using the \emph{Newton
polytope}~\cite{Stu98}, defined as the convex hull of the set of
exponents, considered as vectors in $\R^n$. It was shown by Reznick in
\cite{Rez78} that $Z(x)$ need only contain monomials whose squared
degrees are contained in the convex hull of the degrees of monomials
in $p(x)$.  Consequently, for sparse $p(x)$ the size of the vector
$Z(x)$ and matrix $Q$ appearing in the sum of squares decomposition
can be reduced which results in a decrease of the size of the
semidefinite program.

\begin{figure}
\hspace{\stretch{1}}{\includegraphics[width=0.25\columnwidth]{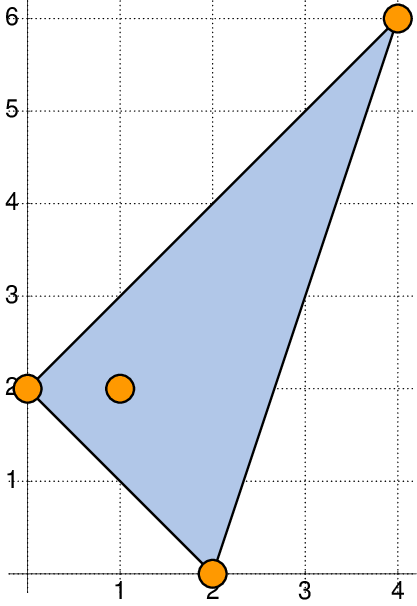}}\hspace{\stretch{1}}
{\includegraphics[width=0.25\columnwidth]{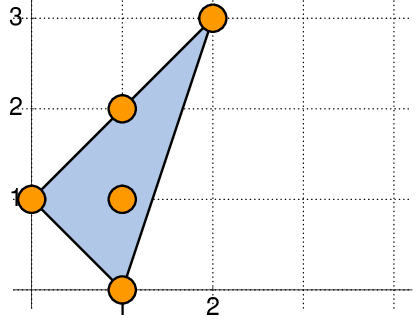}}\hspace{\stretch{1}}
\caption{Newton polytope for the polynomial $p(x,y) = 4 x^4 y^6+x^2-x y^2+y^2$ (left), and possible monomials in its SOS
decomposition (right).}
\label{fig:Newtonpol}
\end{figure}

Consider for example the
polynomial $p(x,y) = 4 x^4 y^6+x^2-x y^2+y^2$, taken from \cite{ParriloLallCDCWorkshop}. Its Newton polytope is
a triangle, being the convex hull of the points
$(4,6),(2,0),(1,2),(2,0)$; see Figure~\ref{fig:Newtonpol}. By the
result mentioned above, we can always find a SOS decomposition that
contains only the monomials $(1,0),(0,1),(1,1),(1,2),(2,3)$.  By
exploiting sparsity, non-negativity of $p(x,y)$ can thus be verified
by solving a semidefinite program of size $5\times 5$ with 13
constraints.  On the other hand, when sparsity is not exploited, we
need to solve a $11\times 11$ semidefinite program with 32
constraints.

When using the argument `sparse', SOSTOOLS takes any sparsity structure into account, and computes an appropriate set of monomials for the sum of squares decomposition to reduce the size of the semidefinite program as described in the above paragraph. To compute the set of monomials defined by the Newton polytope, SOSTOOLS computes the convex hull of a set of monomials either via:
\begin{itemize}
\item the native MATLAB command \texttt{convhulln} (which is based on the software QHULL), or 
\item the specialized external package CDD \cite{CDD}, developed by K. Fukuda. 
\end{itemize}
The choice of the software to be used is determined by the content of the variable \verb"cdd" defined in file \verb"cddpath.m". By default the variable \verb"cdd" is set to be an empty string. This enables the use of \verb"convhulln". To use \verb"cdd", the location of the \verb"cdd" executable file should be assigned to the variable \verb"cdd". Examples are given in commented code within the file \verb"cddpath.m".

%

Special care is taken with the case when the set of
exponents has lower affine dimension than the number of variables
(this case occurs for instance for homogeneous polynomials, where the
sum of the degrees is equal to a constant), in which case a projection
to a lower dimensional space is performed prior to the convex hull
computation.

\subsection{Multipartite Structure}
\label{Sec:Multipartite}

In this section we concentrate on
a particular structure of polynomials that appears frequently in
robust control theory when considering, for instance, Lyapunov
function analysis for linear systems with parametric uncertainty.
For such a case, the indeterminates that
appear in the Lyapunov conditions are Kronecker products of parameters
(zeroth order and higher) and state variables (second order). This
special structure should be taken into account when constructing the
vector $Z(x)$ used in the sum of squares decomposition $p(x) = Z(x)^T
Q Z(x)$. Let us first define what we mean by a \emph{multipartite}
polynomial.

A polynomial $p(x) \in \R[\textbf{x}_1, \ldots,
\textbf{x}_n ]$ in $\sum_{i=1}^n m_i$ indeterminates, where
$\textbf{x}_i = [x_{i1},\ldots,x_{i m_i}]$
given by
\begin{equation*}
p(x) = \sum_{\alpha} c_{\alpha}
\textbf{x}_1^{\alpha_1}\textbf{x}_2^{\alpha_2}\cdots
\textbf{x}_n^{\alpha_n}
\end{equation*}
is termed \emph{multipartite} if for all $i \geq 2$, $
\sum_{k=1}^{m_i} \alpha_{ik}$ is constant, i.e. the monomials in all
but one partition are of \emph{homogeneous} order.
In other words, a multipartite polynomial is homogenous when fixing
any $(n-1)$ blocks of variables, always including the first
block.

This special structure of $p(x)$ can be taken into
account through its Newton polytope.  It has been argued in an earlier
section that when considering the SOS decomposition of a sparse
polynomial (in which many of the coefficients $c_{\alpha}$ are zero),
the nonzero monomials in $Z(x) = [\textbf{x}^{\beta}]$ are the
ones for which $2\beta$ belongs to the convex hull of the degrees
$\alpha$~\cite{Rez78}.  What distinguishes this case from the general
one, is that the Newton polytope of $p(x)$ is the
\emph{Cartesian product} of the individual Newton polytopes
corresponding to the blocks of variables. Hence, the convex hull
should only be computed for the individual $\alpha_i$, which
significantly reduces the complexity and avoids ill-conditioning in
the computation of a degenerate convex hull in a higher dimensional
space.

A specific kind of multipartite polynomials important in practice is
the one that appears when considering \emph{sum of squares
matrices}. These are matrices with polynomial entries that are
positive semi-definite for every value of the indeterminates.
Suppose $S \in \R[\textbf{x}]^{m\times m}$ is a symmetric matrix, and
let  $\textbf{y} = [y_1,\ldots,y_m]$ be new indeterminates. The matrix
$S$ is a \emph{sum of squares (SOS) matrix} if the bipartite
scalar polynomial $\textbf{y}^TS \textbf{y}$ is a sum of squares in
$\R[\textbf{x},\textbf{y}]$.
For example, the matrix $S \in \R[x]^{2 \times 2}$ given by:
\begin{equation*}
S = \left[ \begin{array}{cc}
  x^2-2x+2 & x \\
  x & x^2
\end{array} \right]
\end{equation*}
is a SOS matrix, since
\begin{eqnarray*}
\mathbf{y}^T S \mathbf{y} &=& \left [\begin{array}{c}
  y_1 \\
  x y_1 \\
  y_2 \\
  x y_2
\end{array} \right]^T \left[
\begin{array}{cccc}
  2 & -1 & 0 & 1 \\
  -1 & 1 & 0 & 0 \\
  0 & 0 & 0 & 0 \\
  1 & 0 & 0 & 1
\end{array} \right] \left [\begin{array}{c}
  y_1 \\
  x y_1 \\
  y_2 \\
  x y_2
\end{array} \right] \\ &=& (y_1+xy_2)^2+(xy_1-y_1)^2.
\end{eqnarray*}
Note that SOS matrices for which $n = 1$, i.e. $S \in \R[x]^{m \times
m}$, are positive semidefinite for all real $x$ if and only if they
are SOS matrices; this is because the resulting polynomial will be
second order in $\textbf{y}$ and will only contain one variable $x$;
the resulting positive semidefinite biform is always a sum of
squares~\cite{CLRrealzeros}.

In this manner a SOS matrix in several variables can be converted to a
SOS polynomial, whose decomposition is computed using semidefinite
programming. Because of the bipartite structure, only monomials in the
form $x_i^ky_j$ will appear in the vector $Z$, as mentioned earlier.
For example, the above sum of squares matrix can
be verified as follows:
\begin{matlab}
\begin{verbatim}
>> syms x y1 y2 real;
>> S = [x^2-2*x+2 , x ; x, x^2];
>> y = [y1 ; y2];
>> p = y' * S * y ;
>> prog = sosprogram([x,y1,y2]);
>> prog = sosineq(prog,p,'sparsemultipartite',{[x],[y1,y2]});
>> prog = sossolve(prog);
\end{verbatim}
\end{matlab}

\subsection{Matrix Inequality Constraints}

The function \verb"sosmatrixineq"  allows users to specify polynomial matrix inequalities. Matrix inequality constraints are interpreted in SOSTOOLS as follows.

We recall that for a given symmetric matrix $M\in\R[\theta]^{r\times r}$ we say that $M\succeq 0$ if $y^{T}M(\theta)y\ge 0$ for all $y,\theta$ where $y=[y_{1},\hdots,y_{r}]^{T}.$ 

Alternatively, given a symmetric polynomial matrix $M\in\R[\theta]^{r\times r}$ we say that $M(\theta)$ is a sum of squares matrix if there exist a polynomial matrix $H(\theta) \in \R[\theta]^{s \times r}$ for some $s \in \N$ such that $M(\theta) = H^{T}(\theta)H(\theta)$. This is equivalent to $y^{T}M(\theta)y$ being a sum of squares in $\R[\theta,y]$, \cite{Koj03}.

In SOSTOOLS \emph{both} of the formulations mentioned above are available to the user and are declared using \verb"sosmatrixineq". For a polynomial matrix variable $M(\theta)$, it allows one to either define expressions such as $y^{T}M(\theta)y$ to be a sum of squares polynomial in $\R[y,\theta]$, or to directly constrain the matrix $M(\theta)$ to be a sum of squares  polynomial matrix  in $\R^{r \times r}[\theta]$.

Once a matrix variable $M(\theta)$ has been declared using \verb"sospolymatrixvar", we constrain $y^{T}M(\theta)y$ to be a sum of squares polynomial using \verb"sosmatrixineq" as illustrated below.

\begin{matlab}
\begin{verbatim}
>> syms theta1 theta2;
>> theta = [theta1, theta2];
>> prog = sosprogram(theta);
>> [prog,M] = sospolymatrixvar(prog,monomials(theta,0:2),[3,3],'symmetric');
>> prog = sosmatrixineq(prog,M,'quadraticMineq');
\end{verbatim}
\end{matlab}
Note that the user does \emph{not} declare the vector of indeterminates $y$ as this is handled internally by SOSTOOLS. The naming convention used for the indeterminates is \verb"Mvar_i" where the subscript \verb"i" is used for indexing\footnote{For this reason \texttt{Mvar} should not be used as a variable name.}. In order to save memory, the \verb"Mvar_i" variables may be used with multiple inequalities, therefore \verb"i" will never be greater than the dimension of the largest matrix inequality. 

Some basic error checking has been included that ensures that the matrix argument being passed into \verb"sosmatrixineq" is square and symmetric.

Else, one can ask for $M(\theta)$ to be sum of squares matrix in $\R^{r \times r}[\theta]$ using the following command:

\begin{matlab}
\begin{verbatim}
>> prog = sosmatrixineq(prog,M,'Mineq');
\end{verbatim}
\end{matlab}

A special case of polynomial matrix inequalities are Linear Matrix Inequalities (LMIs). An LMI corresponds to a degree zero polynomial matrix. The following example illustrates how SOSTOOLS can be used as an LMI solver. Consider the LTI dynamical system
\begin{equation}\label{eq:LTI}
\frac{d}{dt}x(t) = Ax(t)
\end{equation}
where $x(t)\in \R^{n}$. It is well known that (\ref{eq:LTI}) is asymptotically stable if and only if there exists a positive definite matrix $P\in\R^{n\times n}$ such that $A^{T}P+PA\prec0$. The code below illustrates how this can be implemented in SOSTOOLS.

\begin{matlab}
\begin{verbatim}
>> syms x;
>> G = rss(4,2,2);
>> A = G.a;
>> eps = 1e-6; I = eye(4);
>> prog = sosprogram(x);
>> [prog,P] = sospolymatrixvar(prog,monomials(x,0),[4,4],'symmetric');
>> prog = sosmatrixineq(prog,P-eps*I,'quadraticMineq');
>> deriv = A'*P+P*A;
>> prog = sosmatrixineq(prog,-deriv-eps*I,'quadraticMineq');
>> prog = sossolve(prog);
>> P = double(sosgetsol(prog,P));
>> A

A =

     -1.2083   -0.6003    0.0488   -0.2103
     -0.6003   -1.6257    0.1917    0.0031
      0.0488    0.1917   -2.1235   -1.0333
     -0.2103    0.0031   -1.0333   -1.6770

>> P

P =

     0.5612   -0.1350    0.0190   -0.0542
    -0.1350    0.4675    0.0288   -0.0003
     0.0190    0.0288    0.4228   -0.1768
    -0.0542   -0.0003   -0.1768    0.4984
\end{verbatim}
\end{matlab}

\textbf{NOTE:} The custom function \verb"findsos" (see example in Section~\ref{sec:SOSTest}) of SOSTOOLS is overloaded to accept matrix arguments. Let  \verb"M" be a given symmetric polynomial matrix of dimension $r\times r$. Calling  \verb"findsos" with the argument \verb"M" will will return the Gram matrix \verb"Q", the vector of monomials \verb"Z" and the decomposition \verb"H" such that
\begin{equation*}
M(x) =H^{T}(x)H(x) = (I_{r}\otimes Z(x))^{T}Q(I_{r}\otimes Z(x)),
\end{equation*}
where $I_{r}$ denotes the $r$-dimensional identity matrix.   In order to illustrate this functionality, consider the following code for which we extract the SOS matrix decomposition from a solution to a SOS program with matrix constraints.

\begin{matlab}
\begin{verbatim}
>> syms x1 x2 x3;
>> x = [x1, x2, x3];
>> prog = sosprogram(x);
>> [prog,M] = sospolymatrixvar(prog,monomials(x,0:2),[3,3],'symmetric');
>> prog = sosmatrixineq(prog,M-(x1^2+x2^2+x3^2)*eye(3),'Mineq');
>> prog = sossolve(prog);
>> M = sosgetsol(prog,M);
>> [Q,Z,H] = findsos(M);
\end{verbatim}
\end{matlab}

\section{Setting an Objective Function}
The function \verb"sossetobj" is used to set an objective function in an optimization problem. The objective function has to be a linear function of the decision variables, and will be minimized by the solver. For instance, if \verb"a" and \verb"b" are symbolic decision variables
in an SOSP named \verb"Program4", then
\begin{matlab}
\begin{verbatim}
>> Program4 = sossetobj(Program4,a-b);
\end{verbatim}
\end{matlab}
will set
\begin{equation}
\text{minimize } (a - b)
\end{equation}
as the objective of \verb"Program4".

Sometimes you may want to minimize an objective function that contains one or more reserved variables
\verb"coeff_nnn", which are created by \verb"sospolyvar" or \verb"sossosvar".
These variables are not individually available in the MATLAB workspace by default. You
must give the argument \verb"'wscoeff'" to the corresponding
\verb"sospolyvar" or \verb"sossosvar" call in order to have these
variables available in the MATLAB workspace. This has been described in Section~2.3.2.

\section{Calling Solver}

A sum of squares program that has been completely defined can be solved using \verb"sossolve.m". If no options are specified, then, for example, to solve \verb"Program5", the command is called with just one argument:
\begin{matlab}
>> Program5 = sossolve(Program5),
\end{matlab}
This function converts the SOSP into an equivalent SDP, calls the default semidefinite programming solver (SeDuMi), and converts the result given by the semidefinite programming solver back into a solution to the original SOSP.  

\subsection{Options} There are several options which can be used when calling \texttt{sossolve}. These options are specified in the \texttt{options} structure and passed to \texttt{sossolve} using the command
\begin{matlab}
>> Program5 = sossolve(Program5,options);
\end{matlab}
Currently, the following options fields are available.

\paragraph{options.solver} The default value for the solver is `SeDuMi'. However, the user may specify other solvers by setting \texttt{options.solver} to one of the following strings.
\begin{itemize}
\item `cdcs'
\item `sdpt3'
\item `csdp'
\item `sdpnal'
\item `sdpnalplus'
\item `sdpa'
\item `mosek'
\end{itemize}
Note that the conversion time to `mosek' input format may be significant for large-scale problems.

\paragraph{options.params} The user may pass a solver-specific parameter structure to any of the supported solvers using the 
 fields of the structure \verb"options.params". For example, when specifying 'SeDuMi' in \verb"options.solver" it is possible to define the tolerance by setting the field \verb"options.params.tol". SeDuMi is called by default with a tolerance of 1e-9. An example illustrating the a user-defined solver and its options is given in Section \ref{sec:EXSOSSetCont}.

\paragraph{options.simplify} Included with SOSTOOLS is the \texttt{sospsimplify} routine, which pre-processes the SDP as described in~\cite{sospsimplify}. For some SOS programming problems, can significantly reduce the size of the resulting SDP. Note that the reductions in \texttt{sospsimplify} are similar to those in `frlib'. To request \texttt{sossolve} to use the \texttt{sospsimplify} routine, the \texttt{options.simplify} field should be specified as follows. 
\begin{matlab}
\begin{verbatim}
>>options.simplify = 'on';
>>prog = sossolve(prog,options);
\end{verbatim}
\end{matlab}

\subsection{Output from \texttt{sossolve}} Typical output that you will get on your screen is shown in Figure~\ref{SOSTOOL Output}. Several things deserve
some explanation:
\begin{itemize}
\item \verb"Size" indicates the size of the resulting SDP.
\item \verb"Residual norm" is the norm of numerical error in the solution.
\item \verb"pinf=1" or \verb"dinf=1" indicate primal or dual infeasibility.
\item \verb"numerr=1" gives a warning of numerical inaccuracy. This is usually accompanied by large
\verb"Residual norm". On the other hand, \verb"numerr=2" is a sign of complete failure because of
numerical problem.
\end{itemize}

\begin{figure}[!htbp]
  \centering
\begin{verbatim}
Size: 10  5

SeDuMi 1.05 by Jos F. Sturm, 1998, 2001.
Alg = 2: xz-corrector, Step-Differentiation, theta = 0.250, beta = 0.500
eqs m = 5, order n = 9, dim = 13, blocks = 3
nnz(A) = 13 + 0, nnz(ADA) = 11, nnz(L) = 8
 it :     b*y       gap    delta  rate   t/tP*  t/tD*   feas cg cg
  0 :            7.00E+000 0.000
  1 : -3.03E+000 1.21E+000 0.000 0.1734 0.9026 0.9000   0.64  1  1
  2 : -4.00E+000 6.36E-003 0.000 0.0052 0.9990 0.9990   0.94  1  1
  3 : -4.00E+000 2.19E-004 0.000 0.0344 0.9900 0.9786   1.00  1  1
  4 : -4.00E+000 1.99E-005 0.234 0.0908 0.9459 0.9450   1.00  1  1
  5 : -4.00E+000 2.37E-006 0.000 0.1194 0.9198 0.9000   0.91  1  2
  6 : -4.00E+000 3.85E-007 0.000 0.1620 0.9095 0.9000   1.00  3  3
  7 : -4.00E+000 6.43E-008 0.000 0.1673 0.9000 0.9034   1.00  4  4
  8 : -4.00E+000 2.96E-009 0.103 0.0460 0.9900 0.9900   1.00  3  4
  9 : -4.00E+000 5.16E-010 0.000 0.1743 0.9025 0.9000   1.00  5  5
 10 : -4.00E+000 1.88E-011 0.327 0.0365 0.9900 0.9905   1.00  5  5
iter seconds digits       c*x               b*y
 10      0.4   Inf -4.0000000000e+000 -4.0000000000e+000
|Ax-b| =  9.2e-011, [Ay-c]_+ =  1.1E-011, |x|= 9.2e+000, |y|= 6.8e+000
Max-norms: ||b||=2, ||c|| = 5,
Cholesky |add|=0, |skip| = 1, ||L.L|| = 2.00001.

Residual norm: 9.2143e-011

       cpusec: 0.3900
         iter: 10
    feasratio: 1.0000
         pinf: 0
         dinf: 0
       numerr: 0
\end{verbatim}
\caption{Output of \name\ (some is generated by SeDuMi).}\label{SOSTOOL Output}
\end{figure}

\section{Getting Solutions}
After your sum of squares program has been solved, you can get the solutions to the program
using \verb"sosgetsol.m". The function takes two arguments, where the first argument is
the SOSP, and the second is a symbolic expression, which typically will be
an SOSP variable. All decision variables in this expression will be substituted by
the numerical values obtained as the solution to the corresponding SDP.
Typing
\begin{matlab}
>> SOLp1 = sosgetsol(Program6,p1);
\end{matlab}
where \verb"p1" is an polynomial variable, for example, will return in \verb"SOLp1"
a polynomial with some numerical coefficients, which is obtained by substituting all decision variables
in \verb"p1" by the numerical solution to the SOSP \verb"Problem6", provided this SOSP has been solved
beforehand.

By default, all the numerical values returned by \verb"sosgetsol" will
have a five-digit presentation.  If needed, this can be changed by
giving the desired number of digits as the third argument to
\verb"sosgetsol", such as
\begin{matlab}
>> SOLp1 = sosgetsol(Program7,p1,12);
\end{matlab}
which will return the numerical solution with twelve digits. Note
however, that this does not change the accuracy of the SDP solution,
but only its presentation.

\chapter{Applications of Sum of Squares Programming}
\label{Applications}

In this chapter we present some problems that can be solved using
\name. The majority of the examples here are from \cite{Par00}, except
when noted otherwise. Many more application examples and customized
files will be included in the near future.
\begin{description}
\item [Note:] For some of the problems here (in particular,
  copositivity and equality-constrained ones such as MAXCUT) the SDP
  formulations obtained by \name\ are not the most efficient ones, as
  the special structure of the resulting polynomials is not fully
  exploited in the current distribution.  This will be incorporated in
  the next release of \name, whose development is already in progress.
\end{description}

\section{Sum of Squares Test}
\label{sec:SOSTest}
As mentioned in Chapter~1, testing if a polynomial $p(x)$ is nonnegative
for all $x\in \mathbb{R}^n$ is a hard problem, but can be relaxed to the
problem of checking if $p(x)$ is an SOS. This can be solved using \name,
by casting it as a feasibility problem.
\myboxed{
\noindent \textbf{SOSDEMO1}:
\vspace{3 mm}\par
\noindent Given a polynomial $p(x)$, determine if
\begin{equation}p(x)\geq 0
\end{equation}
is feasible.}
Notice that even though there are no explicit decision variables in
this SOSP, we still need to solve a semidefinite programming
problem to decide if the program is feasible or not.

The MATLAB code for solving this SOSP can be found in \verb"sosdemo1.m", shown in Figure~\ref{demo1},
and \verb"sosdemo1p.m" (using polynomial objects),
where we consider $p(x)=2x_1^4 + 2x_1^3x_2 - x_1^2x_2^2 + 5x_2^4$. Since the program is feasible,
it follows that $p(x)\geq 0$.

In addition, \name\ provides a function named \verb"findsos" to find an SOS decomposition
of a polynomial $p(x)$. This function returns the coefficient matrix $Q$ and the monomial vector $Z(x)$
which are used in the Gram matrix form. For the same polynomial as above, we may as well type
\begin{matlab}
\begin{verbatim}
>> [Q,Z] = findsos(p);
\end{verbatim}
\end{matlab}
to find $Q$ and $Z(x)$ such that $p(x)=Z^T(x)QZ(x)$. If $p(x)$ is
not a sum of squares, the function will return empty \verb"Q" and
\verb"Z".

For certain applications, it is
particularly important to ensure that the SOS decomposition found
numerically by SDP methods actually corresponds to a true solution,
and is not the result of roundoff errors. This is specially true in
the case of ill-conditioned problems, since SDP solvers can sometimes
produce in this case unreliable results. There are several ways of
doing this, for instance using backwards error analysis, or by
computing rational solutions, that we can fully verify
symbolically. Towards this end, we have incorporated an experimental
option to round to rational numbers a candidate floating point SDP
solution, in such a way to produce an exact SOS representation of the
input polynomial (which should have integer or rational coefficients).
The procedure will succeed if the computed solution is
``well-centered,'' far away from the boundary of the feasible set; the
details of the rounding procedure will be explained elsewhere.

Currently, this facility is available only through the customized
function \texttt{findsos}, by giving an additional input argument
\texttt{`rational'}. On future releases, we may extend this to more
general SOS program formulations.  We illustrate its usage below.
Running
\begin{quote}
\begin{verbatim}
>> syms x y;
>> p = 4*x^4*y^6+x^2-x*y^2+y^2;
>> [Q,Z]=findsos(p,'rational');
\end{verbatim}
\end{quote}
we obtain a rational sum of squares representation for $p(x,y)$ given by {\small
\[
\left[\begin{array}{c} y \\       x \\     x y \\   x
y^2 \\  x^2 y^3
\end{array}\right]^T
\left [\begin {array}{rrrrr}
1&0& -\frac{1}{2}&0&-1\\
0&1&0&-\frac{2}{3}&0\\
-\frac{1}{2}&0&\frac{4}{3}&0&0\\
0&-\frac{2}{3}&0&2&0\\
-1&0&0&0&4
\end {array}\right ]
\left[\begin{array}{c} y \\       x \\     x y \\   x y^2 \\  x^2
y^3
\end{array}\right],
\]}%
where the matrix is given by the symbolic variable
\texttt{Q}, and \texttt{Z} is the vector of
monomials. When polynomial object is used, three output arguments should be given to \texttt{findsos}:
\begin{quote}
\begin{verbatim}
>> pvar x y;
>> p = 4*x^4*y^6+x^2-x*y^2+y^2;
>> [Q,Z,D]=findsos(p,'rational');
\end{verbatim}
\end{quote}
In this case, \verb"Q" is a matrix of integers and \verb"D" is a scalar integer. The variables are related via:
\begin{align*}
p(x,y) = \frac{1}{D}Z^T(x,y)QZ(x,y).
\end{align*}

\begin{figure}[phtb]
\centering
\begin{verbatim}
% SOSDEMO1 --- Sum of Squares Test
% Section 4.1 of SOSTOOLS User's Manual

clear; echo on;
pvar x1 x2;
vartable = [x1, x2];

% =============================================
% First, initialize the sum of squares program
prog = sosprogram(vartable);   % No decision variables.

% =============================================
% Next, define the inequality

% p(x1,x2) >=  0
p = 2*x1^4 + 2*x1^3*x2 - x1^2*x2^2 + 5*x2^4;
prog = sosineq(prog,p);

% =============================================
% And call solver
solver_opt.solver = 'sedumi';
[prog,info] = sossolve(prog,solver_opt);

% =============================================
% If program is feasible, p(x1,x2) is an SOS.
echo off;
\end{verbatim}
\caption{Sum of squares test -- \texttt{sosdemo1.m}}\label{demo1}
\end{figure}

\section{Lyapunov Function Search}
The Lyapunov stability theorem (see e.g. \cite{Kha96}) has been a
cornerstone of nonlinear system analysis for several decades. In
principle, the theorem states that a system $\dot x = f(x)$ with
equilibrium at the origin is stable if there exists a positive
definite function $V(x)$ such that the derivative of $V$ along the
system trajectories is non-positive.

We will now show how to search for Lyapunov function using \name. Consider the system
\begin{equation}
\begin{bmatrix}\dot x_1\\ \dot x_2 \\ \dot x_3\end{bmatrix} =
\begin{bmatrix}-x_1^3 - x_1x_3^2\\ -x_2 -x_1^2x_2\\ -x_3 -\frac{3x_3}{x_3^2+1} + 3x_1^2x_3\end{bmatrix},\label{DynSys}
\end{equation}
with an equilibrium at the origin. Notice that the linearization of
(\ref{DynSys}) has zero eigenvalue, and therefore cannot be used to
analyze local stability of the equilibrium.  Now assume that we are
interested in a quadratic Lyapunov function $V(x)$ for proving
stability of the system. Then $V(x)$ must satisfy
\begin{eqnarray}
V - \epsilon (x_1^2 + x_2^2 + x_3^2) & \geq & 0,\nonumber \\
-\frac{\partial V}{\partial x_1}\dot x_1 -\frac{\partial V}{\partial x_2}\dot x_2
-\frac{\partial V}{\partial x_3}\dot x_3 & \geq & 0. \label{DerivativeEq1}
\end{eqnarray}
The first inequality, with $\epsilon$ being any constant greater than zero, is needed to guarantee
positive definiteness of $V(x)$. However, notice that $\dot x_3$ is a rational function, and therefore
(\ref{DerivativeEq1}) is not a valid SOSP constraint. But since $x_3^2+1 > 0 $ for any $x_3$, we can
just reformulate (\ref{DerivativeEq1}) as
\begin{equation*}
-\frac{\partial V}{\partial x_1}(x_3^2+1)\dot x_1 -\frac{\partial V}{\partial x_2}(x_3^2+1)\dot x_2
-\frac{\partial V}{\partial x_3}(x_3^2+1)\dot x_3 \geq 0.
\end{equation*}
Thus, we have the following SOSP (we choose $\epsilon = 1$):
\myboxed{
\noindent \textbf{SOSDEMO2}:
\vspace{3 mm}\par
\noindent Find a polynomial $V(x)$,
such that
\begin{eqnarray}
V - (x_1^2 + x_2^2 + x_3^2) & \geq & 0,\\
-\frac{\partial V}{\partial x_1}(x_3^2+1)\dot x_1 -\frac{\partial V}{\partial x_2}(x_3^2+1)\dot x_2
-\frac{\partial V}{\partial x_3}(x_3^2+1)\dot x_3 & \geq & 0.
\end{eqnarray}}

The MATLAB code is available in \verb"sosdemo2.m" (or \verb"sosdemo2p.m", when polynomial objects
are used), and is also shown in Figure~\ref{demo2}. The result
given by \name\ is
\begin{equation*}
V(x) = 5.5489x_1^2+4.1068x_2^2+1.7945x_3^2 .
\end{equation*}

The function \verb"findlyap" is provided by \name\ and can be used
to compute a polynomial Lyapunov function for a dynamical system
with polynomial vector field. This function take three arguments,
where the first argument is the vector field of the system, the
second argument is the ordering of the independent variables, and
the third argument is the degree of the Lyapunov function. Thus,
for example, to compute a quadratic Lyapunov function $V(x)$ for
the system
\begin{eqnarray*}
\dot x_1 & = & -x_1^3+x_2,\\
\dot x_2 & = & -x_1-x_2,
\end{eqnarray*}
type
\begin{matlab}
\begin{verbatim}
>> syms x1 x2;
>> V = findlyap([-x1^3+x2; -x1-x2],[x1; x2],2)
\end{verbatim}
\end{matlab}
If no such Lyapunov function exists, the function will return an
empty \verb"V".

\begin{figure}[htbp]
\centering
\begin{verbatim}
% SOSDEMO2 --- Lyapunov Function Search 
% Section 4.2 of SOSTOOLS User's Manual

clear; echo on;
pvar x1 x2 x3;
vars = [x1; x2; x3];

% Constructing the vector field dx/dt = f
f = [(-x1^3-x1*x3^2)*(x3^2+1);
    (-x2-x1^2*x2)*(x3^2+1);
    (-x3+3*x1^2*x3)*(x3^2+1)-3*x3];

% =============================================
% First, initialize the sum of squares program
prog = sosprogram(vars);

% =============================================
% The Lyapunov function V(x): 
[prog,V] = sospolyvar(prog,[x1^2; x2^2; x3^2],'wscoeff');

% =============================================
% Next, define SOSP constraints

% Constraint 1 : V(x) - (x1^2 + x2^2 + x3^2) >= 0
prog = sosineq(prog,V-(x1^2+x2^2+x3^2));

% Constraint 2: -dV/dx*(x3^2+1)*f >= 0
expr = -(diff(V,x1)*f(1)+diff(V,x2)*f(2)+diff(V,x3)*f(3));
prog = sosineq(prog,expr);

% =============================================
% And call solver
solver_opt.solver = 'sedumi';
prog = sossolve(prog,solver_opt);

% =============================================
% Finally, get solution
SOLV = sosgetsol(prog,V)
echo off;
\end{verbatim}
\caption{Lyapunov function search -- \texttt{sosdemo2.m}}\label{demo2}
\end{figure}

\section{Global and Constrained Optimization}
Consider the problem of finding a lower bound for the global minimum
of a function $f(x)$, $x \in \mathbb{R}^n$.  This problem is addressed
in \cite{Sho87}, where an SOS-based approach was first used. A
relaxation method can be formulated as follows. Suppose that there
exists a scalar $\gamma$ such that
\[
f(x)-\gamma \geq 0 \text{ (is an SOS)},
\]
then we know that $f(x) \geq \gamma$, for every $x \in \mathbb{R}^n$.

In this example we will use the Goldstein-Price test function
\cite{GolP71}, which is given by
\begin{eqnarray*}
f(x) & = & [1+(x_1+x_2+1)^2(19-14x_1+3x_1^2-14x_2+6x_1x_2+3x_2^2)]...\\
& & ...\hspace{0.5em}[30+(2x_1-3x_2)^2(18-32x_1+12x_1^2+48x_2-36x_1x_2+27x_2^2)].
\end{eqnarray*}

The SOSP for this problem is
\myboxed{
\vspace{3mm}
\noindent \textbf{SOSDEMO3}:
\vspace{3 mm}\par
\noindent Minimize $-\gamma$,
such that 
\begin{equation}
(f(x)-\gamma) \geq 0.
\end{equation} \vspace{3mm}}

Figure~\ref{demo3} depicts the MATLAB code for this problem. The
optimal value of $\gamma$, as given by \name, is
\[
\gamma_{\text{opt}} = 3.
\]
This is in fact the global minimum of $f(x)$, which is achieved at
$x_1 = 0$, $x_2 = -1$.

\begin{figure}[htbp]
\centering
\begin{verbatim}
% SOSDEMO3 --- Bound on Global Extremum
% Section 4.3 of SOSTOOLS User's Manual 

clear; echo on;
pvar x1 x2 
dpvar gam;
vartable = [x1, x2];

% =============================================
% First, initialize the sum of squares program
prog = sosprogram(vartable);

% =============================================
% Declare decision variable gam too
prog = sosdecvar(prog,gam);

% =============================================
% Next, define SOSP constraints

% Constraint : r(x)*(f(x) - gam) >= 0
% f(x) is the Goldstein-Price function
f1 = x1+x2+1;
f2 = 19-14*x1+3*x1^2-14*x2+6*x1*x2+3*x2^2;
f3 = 2*x1-3*x2;
f4 = 18-32*x1+12*x1^2+48*x2-36*x1*x2+27*x2^2;

f = (1+f1^2*f2)*(30+f3^2*f4);

prog = sosineq(prog,(f-gam));

% =============================================
% Set objective : maximize gam
prog = sossetobj(prog,-gam);

% =============================================
% And call solver
solver_opt.solver = 'sedumi';
prog = sossolve(prog,solver_opt);

% =============================================
% Finally, get solution
SOLgamma = sosgetsol(prog,gam)
echo off
\end{verbatim}
\caption{Bound on global extremum -- \texttt{sosdemo3.m}}\label{demo3}
\end{figure}

The function \verb"findbound" is provided by \name\ and can be used to
find a global lower bound for a polynomial. This function takes just
one argument, the polynomial to be minimized.  The function will
return a lower bound (which may possibly be $-\infty$), a vector
with the variables of the polynomial, and, if an additional condition
is satisfied (the dual solution has rank one), also a point where the
bound is achieved.  Thus, for example, to compute a global minimum for
the polynomial:

\[
F = (a^4+1)(b^4+1)(c^4+1)(d^4+1) +2a + 3 b + 4 c + 5 d,
\]
you would type:
\begin{matlab}
\begin{verbatim}
>> syms a b c d;
>> F = (a^4+1)*(b^4+1)*(c^4+1)*(d^4+1) + 2*a + 3*b + 4*c + 5*d;
>> [bnd,vars,xopt] = findbound(F)
\end{verbatim}
\end{matlab}
For this problem (a polynomial of total degree 16 in
four variables), \name\ returns a certified lower bound
(\verb'bnd=-7.759027')
and also the corresponding optimal point in less than thirty seconds.

In the current version, \texttt{findbound} can also be used
to compute bounds for constrained polynomial optimization problems of the form:
\begin{align*}
& \text{minimize } f(x)\\
& \text{subject to } g_i(x) \geq 0, \quad i = 1,...,M\\
& \hspace{4.6em} h_j(x) = 0, \quad j = 1,...,N.
\end{align*}
A lower bound for $f(x)$ can be computed using
Positivstellensatz-based relaxations.  Assume that there exists a set
of sums of squares $\sigma_j(x)$'s, and a set of polynomials
$\lambda_i(x)$'s, such that
\begin{align}
f(x) - \gamma & =  \sigma_0(x) + \sum_j \lambda_j(x) h_j(x) + \sum_i
\sigma_i(x)g_i(x) + \sum_{i_1,i_2} \sigma_{i_1,i_2}(x)g_{i_1}(x)g_{i_2}(x) + \cdots,
\label{Eq:Positivstellensatz}
\end{align}
then it follows that $\gamma$ is a lower bound for the constrained
optimization problem stated above. This specific kind of
representation corresponds to Schm\"udgen's theorem
\cite{Schmudgen}. By maximizing $\gamma$, we can obtain a lower bound
that becomes increasingly tighter as the degree of the expression
(\ref{Eq:Positivstellensatz}) is increased.

As an example, consider the problem of minimizing $x_1+x_2$,
subject to $x_1\geq 0$, $x_2 \geq 0.5$, $x_1^2+x_2^2 = 1,x_2-x_1^2-0.5=0$. A lower
bound for this problem can be computed using SOSTOOLS as follows:
\begin{quote}
\begin{verbatim}
>> syms x1 x2;
>> degree = 4;
>> [gam,vars,opt] = findbound(x1+x2,[x1, x2-0.5],...
      [x1^2+x2^2-1, x2-x1^2-0.5],degree);
\end{verbatim}
\end{quote}
In the above command, \texttt{degree} is the desired degree for the
expression (\ref{Eq:Positivstellensatz}).  The function
\texttt{findbound} will automatically form the products
$g_{i_1}(x)g_{i_2}(x)$, $g_{i_1}(x)g_{i_2}(x)g_{i_3}(x)$ and so on;
and then construct the sum of squares and polynomial multiplier
$\sigma(x)$'s, $\lambda(x)$'s, such that the degree of the whole
expression is no greater than \texttt{degree}. For this example, a
lower bound of the optimization problem is \texttt{gam}$=1.3911$
corresponding to the optimal solution $x_1 = 0.5682$, $x_2 = 0.8229$,
which can be extracted from the output argument \texttt{opt}.

\section{Matrix Copositivity}
The matrix copositivity problem can be stated as follows:\vspace{2mm}
\par
\begin{minipage}[c]{0.85 \textwidth}
Given a matrix $J\in \mathbb{R}^{n\times n}$,
check if it is copositive, i.e.\ if $y^TJy\geq 0$ for all $y\in \mathbb{R}^n$, $y_i\geq 0$.
\end{minipage}\vspace{2mm} \par \noindent
It is known that checking copositivity of a matrix is a co-NP complete
problem. However, there exist computationally tractable relaxations
for copositivity checking. One relaxation \cite{Par00} is performed by
writing $y_i=x_i^2$, and checking if
\begin{equation}
\left( \sum_{i=1}^{n}x_i^2\right)^m \begin{bmatrix}x_1^2\\ \vdots \\ x_n^2\end{bmatrix}^T
J\begin{bmatrix}x_1^2\\ \vdots \\ x_n^2\end{bmatrix} \triangleq R(x) \label{Rx}
\end{equation}
is an SOS.

Now consider the matrix
\[
J = \left[\begin{array}{rrrrr} 1 & -1 & 1 & 1 & -1\\ -1& 1 & -1 & 1 & 1\\ 1& -1 & 1 & -1 & 1\\ 1 & 1 & -1 & 1 & -1\\
-1 & 1 & 1 & -1 & 1\end{array}\right].
\]
It is known that the matrix above is copositive. This will be proven
using \name.  For this purpose, we have the following SOSP.
\myboxed{
\noindent \textbf{SOSDEMO4}:
\vspace{3 mm}\par
\noindent Determine if
\begin{equation}
R(x) \geq 0,
\end{equation}
is feasible, where $R(x)$ is as in (\ref{Rx}).}

Choosing $m=0$ does not prove that $J$ is copositive. However,
DEMO4 is feasible for $m=1$, and therefore it proves that $J$ is copositive.
The MATLAB code that implements this is given in \verb"sosdemo4.m" and shown
in Figure~\ref{demo4}.

\begin{figure}[htbp]
\centering
\begin{verbatim}
% SOSDEMO4 --- Matrix Copositivity
% Section 4.4 of SOSTOOLS User's Manual

clear; echo on;
pvar x1 x2 x3 x4 x5;
vartable = [x1; x2; x3; x4; x5];

% The matrix under consideration
J = [1 -1  1  1 -1;
    -1  1 -1  1  1;
     1 -1  1 -1  1;
     1  1 -1  1 -1;
    -1  1  1 -1  1];

% =============================================
% First, initialize the sum of squares program
prog = sosprogram(vartable);     % No decision variables.

% =============================================
% Next, define SOSP constraints

% Constraint : r(x)*J(x) - p(x) = 0
J = [x1^2 x2^2 x3^2 x4^2 x5^2]*J*[x1^2; x2^2; x3^2; x4^2; x5^2];
r = x1^2 + x2^2 + x3^2 + x4^2 + x5^2;

prog = sosineq(prog,r*J);

% =============================================
% And call solver
solver_opt.solver = 'sedumi';
prog = sossolve(prog,solver_opt);

% =============================================
% If program is feasible, the matrix J is copositive.
echo off
\end{verbatim}
\caption{Matrix copositivity -- \texttt{sosdemo4.m}}\label{demo4}
\end{figure}

\section{Upper Bound of Structured Singular Value}
Now we will show how \name\ can be used for computing upper bound of
structured singular value $\mu$, a crucial object in
robust control theory (see e.g.\ \cite{DulP00,PacD93}). The following conditions
can be derived from Proposition 8.25 of \cite{DulP00} and
Theorem 6.1 of \cite{Par00}. Given a matrix $M\in\mathbb{C}^{n\times n}$
and structured scalar uncertainties
\[
\Delta = \text{diag}(\delta_1,\delta_2,...,\delta_n),\hspace{2em}\delta_i\in\mathbb{C},
\]
the structured singular value $\mu (M,\Delta)$ is less than $\gamma$, if there exists
solutions $Q_i\geq 0 \in \mathbb{R}^{2n\times 2n},T_i\in \mathbb{R}^{2n\times 2n}$ and $r_{ij}\geq 0$
such that
\begin{equation}
-\sum_{i=1}^{n}Q_i(x)A_i(x) - \sum_{1\leq i < j \leq n}r_{ij}A_i(x)A_j(x) + I(x) \geq 0, \label{muCond}
\end{equation}
where $x\in \mathbb{R}^{2n}$,
\begin{eqnarray}
Q_i(x) & = & x^TQ_ix,\\
I(x) & = & -\sum_{i=1}^{2n}x_i^2,\label{Ix}\\
A_i(x) & = & x^TA_ix,\\
A_i & = & \left[\begin{array}{rr}\text{Re}(H_i) & -\text{Im}(H_i)\\
\text{Im}(H_i) & \text{Re}(H_i)\end{array}\right],\\
H_i & = & M^*e_i^*e_iM - \gamma^2e_i^*e_i,\label{endAix}
\end{eqnarray}
and $e_i$ is the $i$-th unit vector in $\mathbb{C}^n$.

Thus, the SOSP for this problem can be formulated as follows.
\myboxed{
\noindent \textbf{SOSDEMO5}:
\vspace{3 mm}\par
\noindent Choose a fixed value of $\gamma$. For $I(x)$ and $A_i(x)$
as described in (\ref{Ix}) -- (\ref{endAix}), find sums of squares
\begin{eqnarray*}
Q_i(x) & = & x^TQ_ix,\hspace{2em}\text{for }i=1,...,2n,\\
r_{ij} & \geq & 0\hspace{1em}\text{(zero order SOS), for }1\leq i < j \leq 2n,
\end{eqnarray*}
such that (\ref{muCond}) is satisfied.}

The optimal value of $\gamma$ can be found for example by bisection. In \verb"sosdemo5.m"
(Figures~\ref{demo5a}--\ref{demo5b}), we consider the following $M$ (from \cite{PacD93}):
\begin{eqnarray*}
M = UV^*, \hspace{2em}U =\begin{bmatrix} a & 0\\ b & b\\ c & jc\\ d & f\end{bmatrix},
\hspace{2em}V = \begin{bmatrix} 0 & a\\ b & -b\\ c & -jc\\ -jf & -d\end{bmatrix},&&
\end{eqnarray*}
with $a=\sqrt{2/\alpha}$, $b=c=1/\sqrt{\alpha}$, $d=-\sqrt{\beta/\alpha}$, $f = (1+j)\sqrt{1/(\alpha\beta)}$,
$\alpha = 3+\sqrt{3}$, $\beta = \sqrt{3}-1$. It is known that $\mu(M,\Delta)\approx 0.8723$.
Using \verb"sosdemo5.m", we can prove that $\mu(M,\Delta)<0.8724$.

\begin{figure}[htbp]
\centering
\begin{verbatim}
% SOSDEMO5 --- Upper bound for the structured singular value mu
% Section 4.5 of SOSTOOLS User's Manual

clear; echo on;
pvar x1 x2 x3 x4 x5 x6 x7 x8;
vartable = [x1; x2; x3; x4; x5; x6; x7; x8];

% The matrix under consideration
alpha = 3 + sqrt(3);
beta = sqrt(3) - 1;
a = sqrt(2/alpha);
b = 1/sqrt(alpha);
c = b;
d = -sqrt(beta/alpha);
f = (1 + i)*sqrt(1/(alpha*beta));
U = [a 0; b b; c i*c; d f];
V = [0 a; b -b; c -i*c; -i*f -d];
M = U*V';

% Constructing A(x)'s
gam = 0.8724;

Z = monomials(vartable,1);
for i = 1:4
    H = M(i,:)'*M(i,:) - (gam^2)*sparse(i,i,1,4,4,1);
    H = [real(H) -imag(H); imag(H) real(H)];
    A{i} = (Z.')*H*Z;
end

% =============================================
% Initialize the sum of squares program
prog = sosprogram(vartable);

% =============================================
% Define SOSP variables

% -- Q(x)'s -- : sums of squares
% Monomial vector: [x1; ... x8]
for i = 1:4
    [prog,Q{i}] = sossosvar(prog,Z);
end
\end{verbatim}
\caption{Upper bound of structured singular value -- \texttt{sosdemo5.m}, part 1 of 2.}\label{demo5a}
\end{figure}

\begin{figure}[htbp]
\centering
\begin{verbatim}

% -- r's -- : constant sum of squares
Z = monomials(vartable,0);
r = cell(4,4);
for i = 1:4
    for j = (i+1):4
        [prog,r{i,j}] = sossosvar(prog,Z,'wscoeff');
    end
end

% =============================================
% Next, define SOSP constraints

% Constraint : -sum(Qi(x)*Ai(x)) - sum(rij*Ai(x)*Aj(x)) + I(x) >= 0
expr = 0;
% Adding term
for i = 1:4
    expr = expr - A{i}*Q{i};
end
for i = 1:4
    for j = (i+1):4
        expr = expr - A{i}*A{j}*r{i,j};
    end
end
% Constant term: I(x) = -(x1^4 + ... + x8^4)
I = -(x1^4+x2^4+x3^4+x4^4+x5^4+x6^4+x7^4+x8^4);
expr = expr + I;

prog = sosineq(prog,expr);

% =============================================
% And call solver
solver_opt.solver = 'sedumi';
prog = sossolve(prog,solver_opt);

% =============================================
% If program is feasible, thus 0.8724 is an upper bound for mu.
echo off
\end{verbatim}
\caption{Upper bound of structured singular value -- \texttt{sosdemo5.m}, part 2 of 2.}\label{demo5b}
\end{figure}

\section{MAX CUT}
We will next consider the MAX CUT problem. MAX CUT is the problem of partitioning
nodes in a graph into two disjoint sets $V_1$ and $V_2$, such that the weighted number
of nodes that have an endpoint in $V_1$ and the other in $V_2$ is maximized. This can
be formulated as a boolean optimization problem
\[
\max_{x_i \in \{-1,1\}} \frac{1}{2}\sum_{i,j}w_{ij}(1-x_ix_j),
\]
or equivalently as a constrained optimization
\[
\max_{x_i^2 = 1}f(x) \triangleq \max_{x_i^2 = 1} \frac{1}{2}\sum_{i,j}w_{ij}(1-x_ix_j).
\]
Here $w_{ij}$ is the weight of edge connecting nodes $i$ and $j$. For
example we can take $w_{ij}=0$ if nodes $i$ and $j$ are not connected,
and $w_{ij}=1$ if they are connected. If node $i$ belongs to $V_1$,
then $x_i=1$, and conversely $x_i=-1$ if node $i$ is in $V_2$.

A sufficient condition for $\max_{x_i^2 = 1}f(x) \leq\gamma$ is as
follows. Assume that our graph contains $n$ nodes. Given $f(x)$ and
$\gamma$, then $\max_{x_i^2 = 1}f(x) \leq\gamma$ if there exists sum
of squares $p_1(x)$ and polynomials $p_2(x),...,p_{n+1}(x)$ such that
\begin{equation}
p_1(x)(\gamma - f(x)) + \sum_{i=1}^{n}\left(p_{i+1}(x)(x_i^2-1)\right) - (\gamma-f(x))^2 \geq 0. \label{MAXCUTCond}
\end{equation}
This can be proved by a contradiction. Suppose there exists $x\in
\{-1,1\}^n$ such that $f(x)>\gamma$.  Then the first term in
(\ref{MAXCUTCond}) will be negative, the terms under summation will be
zero, and the last term will be negative. Thus we have a
contradiction.


For \verb"sosdemo6.m" (see Figure~\ref{demo6}),
we consider the 5-cycle, i.e., a graph with $5$ nodes and $5$ edges forming
a closed chain. The number of cut is given by
\begin{equation}
f(x) = 2.5 - 0.5x_1x_2 - 0.5x_2x_3 - 0.5x_3x_4 - 0.5x_4x_5 - 0.5x_5x_1.\label{fx}
\end{equation}
Our SOSP is as follows.
\myboxed{
\noindent \textbf{SOSDEMO6}:
\vspace{3 mm}\par
\noindent Choose a fixed value of $\gamma$. For $f(x)$ given in (\ref{fx}), find\vspace{3mm}\\
\hspace*{2em}sum of squares $p_1(x)  = \begin{bmatrix}1 \\ x\end{bmatrix}^TQ\begin{bmatrix}1 \\ x\end{bmatrix}$\\
\hspace*{2em}polynomials $p_{i+1}(x)$ of degree 2, for $i = 1,...,n$\vspace{3mm}\\
\noindent such that (\ref{MAXCUTCond}) is satisfied.}

Using \verb"sosdemo6.m", we can show that $f(x)\leq 4$. Four is indeed the maximum cut for $5$-cycle.

\begin{figure}[htbp]
\centering
\small
\begin{verbatim}
% SOSDEMO6 --- MAX CUT
% Section 4.6 of SOSTOOLS User's Manual

clear; echo on;
pvar x1 x2 x3 x4 x5;
vartable = [x1; x2; x3; x4; x5];

% Number of cuts
f = 2.5 - 0.5*x1*x2 - 0.5*x2*x3 - 0.5*x3*x4 - 0.5*x4*x5 - 0.5*x5*x1;

% Boolean constraints
bc{1} = x1^2 - 1 ;
bc{2} = x2^2 - 1 ;
bc{3} = x3^2 - 1 ;
bc{4} = x4^2 - 1 ;
bc{5} = x5^2 - 1 ;

% =============================================
% First, initialize the sum of squares program
prog = sosprogram(vartable);

% =============================================
% Then define SOSP variables

% -- p1(x) -- : sum of squares
% Monomial vector: 5 independent variables, degree <= 1
Z = monomials(vartable,[0 1]); 
[prog,p{1}] = sossosvar(prog,Z);

% -- p2(x) ... p6(x) : polynomials
% Monomial vector: 5 independent variables, degree <= 2
Z = monomials(vartable,0:2);
for i = 1:5
    [prog,p{1+i}] = sospolyvar(prog,Z);
end

% =============================================
% Next, define SOSP constraints

% Constraint : p1(x)*(gamma - f(x)) +  p2(x)*bc1(x)
%               + ... + p6(x)*bc5(x) - (gamma-f(x))^2 >= 0
gamma = 4;
\end{verbatim}
\caption{MAX CUT -- \texttt{sosdemo6.m}, part 1 of 2.}\label{demo6}
\end{figure}

\begin{figure}[htbp]
\centering
\small
\begin{verbatim}

expr = p{1}*(gamma-f);
for i = 2:6
    expr = expr + p{i}*bc{i-1};
end
expr = expr - (gamma-f)^2;

prog = sosineq(prog,expr);

% =============================================
% And call solver
solver_opt.solver = 'sedumi';
prog = sossolve(prog,solver_opt);

% =============================================
% If program is feasible, 4 is an upper bound for the cut.
echo off
\end{verbatim}
\caption{MAX CUT -- \texttt{sosdemo6.m}, part 2 of 2.}\label{demo6}
\end{figure}

\section{Chebyshev Polynomials}
This example illustrates the \texttt{sosineq} range-specification
option for univariate polynomials (see Section 2.4.2), and is based on a
well-known extremal property of the Chebyshev polynomials. Consider
the optimization problem:

\myboxed{
\noindent \textbf{SOSDEMO7}:
\vspace{3 mm}\par
\noindent Let $p_n(x)$ be a univariate polynomial of degree $n$, with
$\gamma$ being the coefficient of $x^n$. \\

Maximize $\gamma$, subject to: \vspace{3mm}\\
\hspace*{2em} $|p_n(x)| \leq 1, \quad \forall x \in [-1,1]$.}

The absolute value constraint can be easily rewritten using two
inequalities, namely:
\[
\begin{array}{c} 1+ p_n(x) \geq 0 \\
 1- p_n(x) \geq 0 \end{array}
, \quad        \forall x \in [-1,1].
\]
The optimal solution is $\gamma^* = 2^{n-1}$, with
$p_n^*(x) = \arccos(\cos n x)$ being the $n$-th Chebyshev polynomial of
the first kind.

Using \verb"sosdemo7.m" (shown in Figure~\ref{demo7}), the problem can be easily solved for small
values of $n$ (say $n \leq 13$), with SeDuMi aborting with numerical
errors for larger values of $n$. This is due to the ill-conditioning
of the problem (at least, when using the standard monomial basis).

\begin{figure}[htbp]
\centering
\begin{verbatim}
% SOSDEMO7s --- Chebyshev polynomials
% Section 4.7 of SOSTOOLS User's Manual

clear; echo on;
syms x gam;

% Degree of Chebyshev polynomial
ndeg = 8;   

% =============================================
% First, initialize the sum of squares program
prog = sosprogram([x],[gam]);

% Create the polynomial P
Z = monomials(x,[0:ndeg-1]);
[prog,P1] = sospolyvar(prog,Z);
P = P1 + gam * x^ndeg;           % The leading coeff of P is gam

% Imposing the inequalities
prog = sosineq(prog, 1 - P, [-1, 1]);
prog = sosineq(prog, 1 + P, [-1, 1]);

% And setting objective
prog = sossetobj(prog, -gam);

% Then solve the program
solver_opt.solver = 'sedumi';
prog = sossolve(prog,solver_opt);

% =============================================
% Finally, get solution
SOLV = sosgetsol(prog,P)
GAM = sosgetsol(prog,gam)
echo off
\end{verbatim}
\caption{Chebyshev polynomials -- \texttt{sosdemo7s.m}.}\label{demo7}
\end{figure}

\section{Bounds in Probability}
In this example we illustrate how the sums of squares programming
machinery can be used to obtain bounds on the worst-case probability
of an event, given some moment information on the distribution. We
refer the reader to the work of Bertsimas and Popescu
\cite{BertsimasPopescu} for a detailed discussion of the general case,
as well as references to earlier related work.

Consider an unknown arbitrary probability distribution $q(x)$, with
support in $x \in [0,5]$. We know that its mean $\mu$ is equal to 1,
and its standard deviation $\sigma$ is equal to $1/2$. The question
is: what is the worst-case probability, over all feasible
distributions, of a sample having $x\geq4$?

Using the tools in \cite{BertsimasPopescu}, it can be shown that a
bound on (or in this case, the optimal) worst case value can be found
by solving the optimization problem:
\myboxed{
\noindent \textbf{SOSDEMO8}:
\vspace{3 mm}\par
Minimize $a m_0 + b m_1 + c m_2$, subject to
\[\left\{
\begin{array}{ll}
a + b x + c x^2 \geq 0, &\quad \forall x \in [0,5] \\
a + b x + c x^2 \geq 1, &\quad \forall x \in [4,5],
\end{array}
\right.\]
where $m_0=1$, $m_1 = \mu$, and $m_2 = \mu^2 + \sigma^2$.}
\noindent The optimization problem above is clearly an SOSP,
and is implemented in \texttt{sosdemo8.m} (shown in
Figure~\ref{demo8}).

\begin{figure}[htbp]
\centering
\small
\begin{verbatim}
% SOSDEMO8s --- Bounds in Probability
% Section 4.8 of SOSTOOLS User's Manual
 
clear; echo on;
syms x a b c;

% The probability adds up to one.
m0 = 1 ;

% Mean
m1  = 1 ;

% Variance
sig = 1/2 ;

% E(x^2)
m2 = sig^2+m1^2; 

% Support of the random variable
R = [0,5];

% Event whose probability we want to bound
E = [4,5];

% =============================================
% Constructing and solving the SOS program
prog = sosprogram([x],[a,b,c]);

P = a + b*x + c*x^2 ;

% Nonnegative on the support
prog = sosineq(prog,P,R);

% Greater than one on the event
prog = sosineq(prog,P-1,E);

% The bound
bnd =  a * m0 + b * m1 + c * m2 ;

% Objective: minimize the bound
prog = sossetobj(prog, bnd) ;

solver_opt.solver = 'sedumi';
prog = sossolve(prog,solver_opt);

% =============================================
% Get solution
BND = sosgetsol(prog,bnd,16)
PP = sosgetsol(prog,P);
echo off;
\end{verbatim}
\caption{Bounds in
probability -- \texttt{sosdemo8s.m}.}\label{demo8}
\end{figure}

The optimal bound, computed from the optimization problem, is equal to
$1/37$, with the optimal polynomial being $a + b x + c x^2 =
\left(\frac{12 x-11}{37} \right)^2$. The worst case probability
distribution is atomic:
\[ q^*(x) = \frac{36}{37} \, \delta(x-{\textstyle \frac{11}{12}}) +
          \frac{1}{37} \, \delta(x-4).
\]
All these values (actually, their floating point approximations) can
be obtained from the numerical solution obtained using \name.
\clearpage
\section{SOS Matrix Decomposition}
This example illustrates how SOSTOOLS v3.00 can be used to determine if an $r\times r$ polynomial matrix $P$ is an SOS matrix. Furthermore, if $P$ is determined to be SOS then it is shown how the matrix decomposition $P(x)=H^{T}(x)H(x)$ can be computed.
\myboxed{
\noindent \textbf{SOSDEMO9}:
\vspace{3 mm}\par
Given a symmetric polynomial matrix $P\in \R[x]^{r\times r}$ determine if $P$ is an SOS matrix and if so compute the polynomial matrix $H(x)$ such that
\begin{equation}
P(x) = H^{T}(x)H(x).
\end{equation}}
The above feasibility problem is implemented in \verb"sosdemo9.m" for the matrix
\begin{equation*}
P(x)=\left[ \begin{array}{cc}  x_1^4 + x_1^2x_2^2 + x_1^2x_3^2 & x_1x_2x_3^2 - x_1^3x_2 - x_1x_2(x_2^2 + 2x_3^2)\\
                                             x_1x_2x_3^2 - x_1^3x_2 - x_1x_2(x_2^2 + 2x_3^2) & x_1^2x_2^2 + x_2^2x_3^2 + (x_2^2 + 2x_3^2)^2 \end{array}\right].
\end{equation*}

The code in Figure 4.10 can be used to compute a matrix decomposition. The \texttt{findsos} function returns the arguments \texttt{Q}, \texttt{Z} and \texttt{Hsol} such that
\begin{equation*}
H(x) = (I_{r}\otimes Z(x))^{T}Q(I_{r}\otimes Z(x))
\end{equation*}
where $I_{r}$ is the $r\times r$ identity matrix, $Q$ is a positive semidefinite matrix and $Z(x)$ is a vector of monomials.
\begin{figure}[htbp]
\centering
\small
\begin{verbatim}
% SOSDEMO9 --- Matrix SOS decomposition
% Section 4.9 of SOSTOOLS User's Manual

clear; echo on;
pvar x1 x2 x3;
vartable = [x1, x2, x3];

% =============================================
% First, initialize the sum of squares program
prog = sosprogram(vartable);   % No decision variables.

% =============================================
% Consider the following candidate sum of squares matrix P(x)
P = [x1^4+x1^2*x2^2+x1^2*x3^2 x1*x2*x3^2-x1^3*x2-x1*x2*(x2^2+2*x3^2); 
     x1*x2*x3^2-x1^3*x2-x1*x2*(x2^2+2*x3^2) x1^2*x2^2+x2^2*x3^2+(x2^2+2*x3^2)^2];  

% Test if P(x1,x2,x3) is an SOS matrix and return H so that P = H.'*H
solver_opt.solver = 'sedumi';
[Q,Z,H] = findsos(P,[],solver_opt);

% Verify that P - H'*H = 0 and P- kron(I,Z)'*Q*kron(I,Z)= 0 to within 
% numerical tolerance.
tol = 1e-8;
cleanpoly(P - H'*H ,tol)
cleanpoly(P- kron(eye(2),Z)'*Q*kron(eye(2),Z), tol)

% Verify that Q >=0
min(eig(Q))

% =============================================
% If program is feasible, P(x1,x2,x3) is an SOS matrix.
echo off;
\end{verbatim}
\caption{Computing a SOS matrix decomposition -- \texttt{sosdemo9.m}.}\label{demo9}
\end{figure}

\newpage
\section{Set Containment}
\label{sec:EXSOSSetCont}
This example illustrates how SOSTOOLS v3.01 can be used to
compute the entries of a   polynomial matrix $P$ such that it is an
SOS matrix.

It has been shown in \cite{ValTG13} that if the matrix $P(x)$ given by
\begin{equation}
P(x) = \left[\begin{array}{cc}  \theta^2-s(x)(\gamma-p(x)) &
g_0(x)+g_1(x) \\ g_0(x)+g_1(x) & 1 \end{array}\right],
\label{eq:setcont}
\end{equation}
is an SOS matrix, then the following set containment holds:
\begin{equation}
\left\lbrace  x \in \R^2 | p(x) \leq \gamma \right\rbrace \subseteq
\left\lbrace  x \in \R^2 | ( (g_0(x)+g_1(x)) + \theta )(\theta -
(g_0(x)+g_1(x)) ) \geq 0 \right\rbrace.
\end{equation}
where given are $p(x)$, a positive polynomial, $g_{0}\in \R[x]$, and $\theta,\gamma>0$ are positive scalars. If a polynomial $g_1(x)$ and SOS multiplier $s(x)$ are found, then the set containment holds. This problem is a sum of squares feasability problem, the code for this demo is given in Figure 4.11.

\begin{figure}[htbp]
\myboxed{
\noindent \textbf{SOSDEMO10}:
\vspace{3 mm}\par
Given $p(x)\in \R[x]$, $g_0(x)\in \R[x]$, $\theta\in \R$, $\gamma \in
\R$, find \vspace{3mm}\\
\hspace*{2em}polynomial $g_1(x)\in \R[x]$\\
\hspace*{2em}Sum of Squares $s(x)$\vspace{3mm}\\
\noindent such that \eqref{eq:setcont} is a sum of squares matrix.
}
\end{figure}

\begin{figure}[htbp]
\centering
\small
\begin{verbatim}
% SOSDEMO10 --- Set containment
% Section 4.10 of SOSTOOLS User's Manual

clear; echo on;
pvar x1 x2; 
vartable = [x1 x2];

eps = 1e-6;

% =============================================
% This is the problem data
p = x1^2+x2^2;
gamma = 1;
g0 = [2 0]*[x1;x2];
theta = 1;

% =============================================
% Initialize the sum of squares program
prog = sosprogram(vartable);

% =============================================
% The multiplier
Zmon = monomials(vartable,0:4);
[prog,s] = sospolymatrixvar(prog,Zmon,[1 1]);

% =============================================
% Term to be added to g0
Zmon = monomials(vartable,2:3);
[prog,g1] = sospolymatrixvar(prog,Zmon,[1 1]);

% =============================================
% The expression to satisfy the set containment
Sc = [theta^2-s*(gamma-p) g0+g1; g0+g1 1];

prog = sosmatrixineq(prog,Sc-eps*eye(2));

solver_opt.solver = 'sedumi';
prog = sossolve(prog,solver_opt);

s = sosgetsol(prog,s);
g1 = sosgetsol(prog,g1);

% Display function g1 after removing small coefficients
cleanpoly(g1,1e-4)

% =============================================
% If program is feasible, { x |((g0+g1) + theta)(theta - (g0+g1)) >=0 } contains { x | p <= gamma }
echo off;
\end{verbatim}
\caption{Set containment -- \texttt{sosdemo10.m}.}\label{demo10}
\end{figure}

The feasibility test above is formulated and solved in
\verb"sosdemo10.m"  for $p(x) = x_1^2+
x_2^2$, $\gamma = \theta = 1$ and $g_0 = 2x_1$, a sum of squares
variable $s(x)$ of degree $4$ and a polynomial variable $g_1(x)$
containing monomials of degrees 2 and 3. This example illustrates the
use of function \verb"sosineq" having a matrix as an input argument.

\chapter{Interfaces to Additional Packages}\label{ch:interface}
The following packages have been written by researchers to extend the functionality of SOSTOOLS. Below is a brief description of these packages and instructions on how to install them. Please refer to their official documentation for the full details. For technical queries please contact the relevant authors. 
\section{INTSOSTOOLS}

The INTSOSTOOLS package (available from \verb"https://github.com/gvalmorbida/INTSOSTOOLS")  is a plug-in for SOSTOOLS for the formulation of optimization problems subject to one-dimensional integral inequalities such as
\begin{equation}
\label{eq:optproblem}
\begin{array}{l}\mbox{maximize} \ \lambda \\ \mbox{subject to}  \  \int_{0}^{1} \left(f(\theta,u(\theta)) -\lambda \right) \mathrm{d}\theta \geq 0.
 \end{array}
\end{equation}

For polynomial problem data, these optimization problems can be solved using semidefinite programming via SOSTOOLS. The functionalities of the package and examples  are detailed in \cite{VP15}.

\section{frlib}
The frlib package \cite{PerP14} (available from \verb"https://github.com/frankpermenter/frlib") provides a pre-processing step that performs a facial reduction on the positive semidefinite cone in order to produce a simplified SDP. Typically finding the appropriate face (that contains the feasible set) of the cone to optimize over is itself an SDP, moreover it is often too expensive to actually solve this SDP as a pre-processing step. However frlib provides a method that optimizes over a specific approximation of the cone that leads to a subset of faces that can be optimized over. Currently the approximations available are \emph{non-negative diagonal} $(\cD^n)$ and diagonally dominant $(\cD\cD^n) $ approximations, where
\begin{equation*}
\cD^n := \left\{ Q\in \mathbb S^n ~|~ Q_{ii}\ge 0,~ Q_{ij}=0 \text{ if }i\neq j \right\}
\end{equation*}
and 
\begin{equation*}
\cD\cD^n := \left\{ Q\in \mathbb S^n ~|~ Q_{ii} \ge \sum_{j\neq i}^n |Q_{ij}| \right\}.
\end{equation*}
approximate the cone of $n\times n$ positive semidefinite matrices. It can be shown that $\cD^n \subset \cD\cD^n$.
\subsection{Example code}
In order to invoke a call to frlib the \verb"options" command must be used as illustrated in the code segment below. The user must choose whether to use the $\cD^n$ or $\cD\cD^n$ approximation of the positive semidefinite cone using the code \verb"options.frlib.approx = x" where \verb"x" is either \verb"'d'" or \verb"'dd'". Additionally it can be specified whether or not to use a QR decomposition in the numerical implementation. The code below selects the $\cD\cD^n$ approximation, a QR decomposition, and sets the solver as SeDuMi: 

\begin{matlab}
\begin{verbatim}
>>options.frlib.approx = 'dd';
>>options.frlib.useQR = 1;
>>options.solver = 'sedumi';
>>prog = sossolve(prog,options);
\end{verbatim}
\end{matlab}

The output shown below indicates that frlib has found a reduction and the solver has been called.

\begin{figure}[!h]
  \centering
\begin{verbatim}
Using frlib toolbox for additional pre-processing...
-------------------------------------------------------------------------
frlib: reductions found!
-------------------------------------------------------------------------
  Dim PSD constraint(s) (original):  3 13 
  Dim PSD constraint(s) (reduced):   3 10 
  
  SeDuMi 1.32 by AdvOL, 2005-2008 and Jos F. Sturm, 1998-2003.
Alg = 2: xz-corrector, Adaptive Step-Differentiation, theta = 0.250, beta = 0.500
Put 3 free variables in a quadratic cone
eqs m = 44, order n = 16, dim = 114, blocks = 4
nnz(A) = 76 + 0, nnz(ADA) = 1540, nnz(L) = 792
 it :     b*y       gap    delta  rate   t/tP*  t/tD*   feas cg cg  prec
    0 :            1.92E+00 0.000
  1 :   2.94E+00 3.97E-01 0.000 0.2072 0.9000 0.9000  -0.87  1  1  3.1E+00
  2 :   1.84E+00 7.36E-02 0.000 0.1856 0.9000 0.9000   0.10  1  1  9.1E-01
  3 :   4.99E-02 1.48E-03 0.000 0.0200 0.9900 0.9900   0.78  1  1  2.0E-02
  4 :   4.15E-04 1.23E-05 0.000 0.0083 0.9990 0.9990   1.00  1  1  1.7E-04
  5 :   3.45E-06 1.02E-07 0.000 0.0083 0.9990 0.9990   1.00  1  1  1.4E-06
  6 :   2.86E-08 3.71E-10 0.000 0.0036 0.9990 0.9990   1.00  1  1  6.4E-09

iter seconds digits       c*x               b*y
  6      0.2   Inf  0.0000000000e+00  2.8644510334e-08
|Ax-b| =   2.0e-09, [Ay-c]_+ =   1.1E-09, |x|=  4.2e+01, |y|=  2.6e-08

Detailed timing (sec)
   Pre          IPM          Post
1.365E-01    1.969E-01    4.535E-02    
Max-norms: ||b||=1, ||c|| = 0,
Cholesky |add|=0, |skip| = 0, ||L.L|| = 7.95584.
---------
Facial reduction applied using frlib: Primal and Dual solution available
---------

\end{verbatim}
\caption{SOSTOOLS output when using frlib - \texttt{sosdemo2.m}.}
\end{figure}

\chapter{Inside SOSTOOLS}
In this chapter the data structures that underpin SOSTOOLS are described. The information in this section can help advanced users to manipulate the sosprogram structure in order to create additional functionality. It is assumed from this point on that the user has a strong working knowledge of the functions described in the previous chapter.

As described in Chapter \ref{ch:SOSP} an SOSP is initialized using the command
\begin{matlab}
\begin{verbatim}
>> syms x y z;
>> prog = sosprogram([x;y;z]);
\end{verbatim}
\end{matlab}
The command above will initialize an empty SOSP called \verb"prog" with variables $x,y$ and $z$ and returns the following structure:
\begin{matlab}
\begin{verbatim}
>> prog

prog =

           var:  [1x1 struct]
          expr:  [1x1 struct]
      extravar:  [1x1 struct]
     objective:  []
       solinfo:  [1x1 struct]
      vartable:  '[x,y,z]'
   symvartable:  [3x1 sym]
        varmat:  []
   decvartable:  '[]'
symdecvartable:  []
\end{verbatim}
\end{matlab}

The various fields above are populated by the addition of new decision variables (polynomials and sum of squares polynomials), constraints (equality and inequality\footnote{Recall that the inequality $h(x)\ge 0$ is interpreted as $h(x)$ having a sum of squares decomposition.}) and the inclusion of an objective function. The contents and structure of each of these fields will now be described and specifically related to the construction of the underlying SDP.

We will illustrate the program structure by constructing a typical sum of squares program that is similar to \verb"SOSDEMO2". The first fields to be populated are \verb"prog.symvartable" and \verb"prog.vartable".\footnote{It is assumed that the symbolic toolbox is being used.} The following output will be seen:
\begin{matlab}
\begin{verbatim}
>> prog.symvartable

prog.symvartable =
    x
    y
    z

 >> prog.vartable

prog.vartable =
    [x,y,z]
\end{verbatim}
\end{matlab}
Note that if the symbolic toolbox is not being used then the \verb"prog.symvartable" field will not exist and \verb"prog.vartable" (which is a character array) is used instead.

Next, a polynomial variable
\verb"V" in the monomials $x^{2},y^{2},z^{2}$ is declared using
\begin{matlab}
\begin{verbatim}
>> [prog,V] = sospolyvar(prog,[x^2; y^2; z^2]);
\end{verbatim}
\end{matlab}
which sets the field \verb"prog.var" as follows:
\begin{matlab}
\begin{verbatim}
>> prog.var

prog.var =

     num:   1
    type:   {'poly'}
       Z:   {[3x3 double]}
      ZZ:   {[3x3 double]}
       T:   {[3x3 double]}
     idx:   {[1]  [4]}
\end{verbatim}
\end{matlab}
The field \verb"prog.var.num" is an integer that gives the number of variables declared. In this case there is only one variable \verb"V", as such each of the remaining fields (excluding \verb"prog.var.idx") contains a single entry. As more variables are declared these fields will contain arrays where each element corresponds to exactly one variable. The type of variable, i.e. polynomial or sum of squares polynomial, is indicated by \verb"prog.var.type": for example an SOSP with two polynomials and a sum of squares polynomial variable would result in
\begin{matlab}
\begin{verbatim}
>> prog.var.type

prog.var.type =

     'poly'  'poly'   'sos'
\end{verbatim}
\end{matlab}
Returning to the example, recall that the polynomial \verb"V" consists of three monomials $x^{2},y^{2},z^{2}$ with unknown coefficients. The field \verb"prog.var.Z" is the monomial degree matrix whose entries are the monomial exponents. For a polynomial in $n$ variables containing $m$ monomials this matrix will have dimension $m\times n$. In this case \verb"V" is clearly a $3\times 3$ matrix with $2$'s on the diagonal. A more illuminating example is given below.
\begin{matlab}
\begin{verbatim}
>> [prog,h] = sospolyvar(prog,[x^2; y^2; z^2; x*y*z; x^2*y]);
>> full(prog.var.Z{2})

ans =
             2  0  0
             0  2  0
             0  0  2
             1  1  1
             2  1  0
\end{verbatim}
\end{matlab}
Note that by default this matrix is stored as a sparse matrix. Define the vector of monomials (without coefficients) that describes \verb"V" by $Z$, i.e. $Z=[x^{2},y^{2},z^{2}]^{T}$. Further, define $W$ to be the vector of pairwise different monomials of the entries of the matrix $ZZ^{T}$. The exponents of all such monomials are then stored in the degree matrix \verb"prog.var.ZZ" where the rows corresponding to the unique set of monomials are ordered lexicographically. The field \verb"prog.var.idx" contains indexing information that is required for constructing the SDP. We will expand upon this indexing further when describing the \verb"prog.extravar" fields.

The next field of interest is \verb"prog.expr" which is the primary field where SOSTOOLS stores the user's constraints. Continuing on with the example we see that there are two constraints in the programme and that both are sum of squares. This can be seen by inspecting \verb"prog.expr.num" and \verb"prog.expr.type" respectively.
\begin{matlab}
\begin{verbatim}
>> prog.expr

prog.expr =

     num:   2
    type:   {'ineq', 'ineq'}
      At:   {[3x3 double]   [3x10 double]}
       b:   {[3x1 double]   [10x1 double]}
       Z:   {[3x3 double]   [10x3 double]}
\end{verbatim}
\end{matlab}
Recall that the canonical primal SDP takes the form:
\begin{eqnarray}
\underset{x}{\text{minimize}}&&c^{T}x  \nonumber\\
\text{s.t.}&&Ax = b \label{eq:SDPprimal}\\
&&x\in \cK \nonumber
\end{eqnarray}
where $\cK$ denotes the symmetric cone of positive semidefinite matrices. Here the two inequalities imposed in the SOSP are converted into their SDP primal form where the field \verb"prog.expr.At" refers to the transposition of the matrix $A$ in (\ref{eq:SDPprimal}) and likewise \verb"prog.expr.b" refers to the vector $b$ in (\ref{eq:SDPprimal}). Finally the field \verb"prog.expr.Z" contains the matrices of monomial exponents corresponding to the two inequalities. These take exactly the same form as described above for the field \verb"prog.var.Z".

Thus far in the example we have not described the decision variables, i.e. the unknown polynomial coefficients. This will now be illustrated through the constraint  $h(y,y,z)= V-(x^{2}+y^{2}+z^{2})\ge0$ where $V$ contains only quadratic terms. This inequality is interpreted as a sum of squares inequality of the form
\begin{equation*}
h(x,y,z) = \left[\begin{array}{c}x\\y\\z\end{array}\right]^{T}Q\left[\begin{array}{c}x\\y\\z\end{array}\right],
\end{equation*}
where the decision variable is the positive semidefinite matrix $Q$ consisting of the coefficients of the polynomial $V-(x^{2}+y^{2}+z^{2})$. Here the matrix $Q$ is constrained to be of the form
\begin{equation}\label{eq:egQ}
Q = \left[\begin{array}{ccc} \verb"coeff_1" & 0 & 0\\ 0 & \verb"coeff_2" & 0\\ 0 & 0 & \verb"coeff_3"\end{array}\right]
\end{equation}
and thus the decision variables are the three non-zero coefficients in $Q$. The decision variables can be seen through the \verb"prog.symdecvartable" field
\begin{matlab}
\begin{verbatim}
>> prog.symdecvartable

prog.symdecvartable =

     coeff_1
     coeff_2
     coeff_3
\end{verbatim}
\end{matlab}
In a similar manner \verb"prog.decvartable" contains the same information but stored as a character array.

This particular example is an SOS feasibility problem (i.e. no objective function has been set). An objective function to be minimised can be set using the function \verb"sossetobj" in which case the field \verb"prog.objective" would contain the relevant data. Recall that for both an SDP and SOSP the objective function must be a linear function. SOSTOOLS will automatically set the weighting vector $c$ in (\ref{eq:SDPprimal}), this is exactly what is contained in \verb"prog.objective", while the $x$ corresponds to the decision variables, i.e. the unknown polynomial coefficients in (\ref{eq:egQ}).

Many SOSPs may include matrix inequality constraints, that is constraints of the form\linebreak $x^{T}M(\theta)x\ge 0$, or more accurately that the previous expression is a sum of squares polynomial in $x$ and $\theta$. When setting such constraints the user does not need to declare the independent variable vector $x$ as this is handled internally by SOSTOOLS. The following code sets the matrix constraint:
\begin{matlab}
\begin{verbatim}
>> sym theta
>> [prog,M] = sospolymatrixvar(prog,monomials([theta],0:2),[3 3],'symmetric');
>> prog = sosmatrixineq(prog,M,'quadraticMineq');
\end{verbatim}
\end{matlab}
which then populates the field \verb"prog.varmat".

\begin{matlab}
\begin{verbatim}
>> prog.varmat

 prog.varmat =

       vartable: '[Mvar_1,Mvar_2,Mvar_3]'
    symvartable: [3x1 sym]
          count: 3

\end{verbatim}
\end{matlab}
Here \verb"prog.varmat.symvartable" is the $3\times 1$ vector of symbolic variables \verb"Mvar_1, Mvar_2, Mvar_3" which correspond to the independent variables $x$ in the above example. The field \linebreak \verb"prog.varmat.symvartable" likewise contains a character array of the vector of variables, while \verb"prog.varmat.count" indicates the number of variables used. Note that in order to save memory the variables \verb"Mvar_1" may be used with multiple inequalities.

There is one further field, \verb"prog.extravar" which SOSTOOLS creates. This field has the following entries:
\begin{matlab}
\begin{verbatim}
>> prog.extravar

prog.extravar =

     num:  2
       Z:  {[3x3 double]   [13x3 double]}
      ZZ:  {[6x3 double]   [52x3 double]}
       T:  {[9x6 double]   [169x52 double]}
     idx:  {[4]  [13]  [182]}
\end{verbatim}
\end{matlab}

At this point all the variables have been defined and both constraints and an objective function (if desired) have been set the SOSP is ready to be solved. This is done using the function
\begin{matlab}
\begin{verbatim}
>> prog = sossolve(prog);
\end{verbatim}
\end{matlab}
which calls the SDP solver and then converts the solution data back into SOSP form. However, the original SDP problem and solution data is stored by SOSTOOLS in the field \verb"prog.solinfo" which contains the subfields:
\begin{matlab}
\begin{verbatim}
>> prog.solinfo

prog.solinfo =

              x:  [181x1 double]
              y:  [58x1 double]
            RRx:  [181x1 double]
            RRy:  [181x1 double]
           info:  [1x1 struct]
  solverOptions:  [1x1 struct]    
            var:  [1x1 struct]
       extravar:  [1x1 struct]
         decvar:  [1x1 struct]
\end{verbatim}
\end{matlab}
It is worth remembering that in general users do not need to know how to access, use or even interpret this data as all the polynomial variables can be accessed via the \verb"sosgetsol" function. For example:
\begin{matlab}
\begin{verbatim}
>> V = sosgetsol(prog,V)

V =

6.6582*x^2+4.5965y^2+2.0747*z^2

\end{verbatim}
\end{matlab}
However there may be occasions when the user would like the specific sum of squares decomposition, or indeed the SDP form of the solution. Before describing how this is done, we explain what each field in \verb"prog.solinfo" contains. Recall that the dual canonical form of an SDP takes the form
\begin{eqnarray}
\underset{y}{\text{maximize}}&&b^{T}y  \nonumber\\
\text{s.t.}&&c-A^{T}y \in \cK \label{eq:SDPdual}.
\end{eqnarray}
Intuitively \verb"prog.solinfo.x" and \verb"prog.solinfo.y" contain the primal and dual decision variables $x$ and $y$ respectively. Note, however, that \verb"prog.solinfo.x" and \verb"prog.solinfo.y" contain the solution vectors to the whole SDP which will typically be a concatenation of multiple smaller SDPs for each constraint. For the SOSDEMO2 example the field \verb"prog.solinfo.x" will contain not only the coefficients to the polynomial $V$ but also the coefficients of the positive semidefinite matrix corresponding to the negativity of the derivative condition. The matrix coefficients are stored in vectorized form. To extract the coefficients of $V$ directly one can use the field \verb"prog.solinfo.var.primal":
\begin{matlab}
\begin{verbatim}
>> prog.solinfo.var.primal{:}

ans =

   6.6582
   4.5965
   2.0747
\end{verbatim}
\end{matlab}
Clearly these are verified as the coefficients of $V$ as given above. For the case where there are multiple variables \verb"prog.solinfo.var.primal" will return an array of vectors. The dual variables, $y$ from (\ref{eq:SDPdual}), can be accessed in much the same way using \verb"prog.solinfo.var.dual". Indeed it is the field \verb"prog.var.idx" and \verb"prog.extravar.idx" that were alluded to earlier that extract the relevant coefficients for each of the polynomial decision variables.

In certain instances it may be desirable to obtain the sum of squares decomposition in order to certify a solution. For example, one may wish to see the sum of squares decomposition, the vector $Z_{i}(x)$ and positive semidefinite matrix $Q_{i}$, of the constraint
\begin{equation*}
-\frac{\partial V}{\partial x_1}(x_3^2+1)\dot x_1 -\frac{\partial V}{\partial x_2}(x_3^2+1)\dot x_2
-\frac{\partial V}{\partial x_3}(x_3^2+1)\dot x_3  \geq  0.
\end{equation*}
This is achieved using the \verb"prog.solinfo.extravar.primal" field. The \verb"prog.solinfo.extravar" field has the following structure:
\begin{matlab}
\begin{verbatim}
>> prog.solinfo.extravar

ans =

    primal:   {[3x3 double]   [13x13 double]}
    double:   {[3x3 double]   [13x13 double]}
\end{verbatim}
\end{matlab}

The derivative condition was the second constraint to be set so in order to obtain the corresponding $Q_{i}$ the following command is used:
\begin{matlab}
\begin{verbatim}
>> Q2 = prog.solinfo.extravar.primal{2};
\end{verbatim}
\end{matlab}
In this example \verb"Q2" is a $13 \times 13$ matrix which we will omit due to space constraints. The matrix $Q_{1}$ corresponding to the first constraint is:
\begin{matlab}
\begin{verbatim}
>> Q1 = prog.solinfo.extravar.primal{1}

Q1 =

   5.6852  0.0000  0.0000
   0.0000  3.5965  0.0000
   0.0000  0.0000  1.0747
\end{verbatim}
\end{matlab}
and the corresponding vector of monomials, $Z_{1}(x,y,z)$ is obtained using
\begin{matlab}
\begin{verbatim}
>> Z1 = mysympower([x,y,z],prog.extravar.Z{1})

Z1 =

   x
   y
   z
\end{verbatim}
\end{matlab}
which proves that $Z_{1}^{T}(x,y,z)Q_{1}Z_{1}(x,y,z) = V(x,y,z)-(x^{2}+y^{2}+z^{2})\ge 0$.

In this example, as we mentioned earlier, there is no objectve function to minimise. However, for SOSPs that do contain an objective function the field \verb"prog.solinfo.decvar" returns both the primal solution $c^{T}x$ corresponding to (\ref{eq:SDPprimal}) and the dual solution $b^{T}y$ corresponding to (\ref{eq:SDPdual}).

The numerical information regarding the solution as returned from the SDP solver is stored in the field \verb"prog.solinfo.info". The exact information stored is dependent upon which of the solvers is used. However the typical information returned is displayed in Figure 3.2.

The fields \verb"prog.solinfo.RRx" and \verb"prog.solinfo.RRy" are populated with the complete solution of the SDP.  The elements of the fields \verb"prog.solinfo.decvar.primal", \linebreak \verb"prog.solinfo.extravar.primal", and\verb"prog.solinfo.decvar.dual",  \verb"prog.solinfo.extravar.dual", are actually computed by extracting and rearranging the entries of \verb"prog.solinfo.RRx" and\linebreak \verb"prog.solinfo.RRy".

The used solver and its parameters are stored in the fields \verb"prog.solinfo.solverOptions.solver" and \verb"prog.solinfo.solverOptions.params".

The notation $RRx$ and $RRy$ is due to the fact that they receive the values of $RR x$, $RR (c - A^{T}y)$ where $x$  and $y$ are the primal and the dual solutions computed by the solver. The matrix $RR$ is a permutation matrix used to rearrange the input data for the solver. This is required because the solvers need to distinguish between decision variables which are on the cone of positive semi-definite matrices and other decision variables. With the permutation matrices $RR$, this arrangement of decision variables previous to the call to the SDP solver is transparent to the user of SOSTOOLS.

\chapter{List of Functions}
{\tt

\begin{tabular}{ll}
\% &  \hspace{-15pt} SOSTOOLS --- Sum of Squares Toolbox \\
\% &  \hspace{-15pt} Version 4.00, 14 September 2021.  \\
\% &\\
\%  &  \hspace{-15pt} Monomial vectors construction:
\end{tabular}

\hspace{-20pt}\begin{tabular}{llcl}
\% &MONOMIALS & --- &Construct a vector of monomials with prespecified \\
\% & & & degrees.\\
\%& MPMONOMIALS & --- & Construct a vector of multipartite monomials with  \\
\% & & & prespecified degrees.
\end{tabular}

\hspace{-20pt}\begin{tabular}{ll}
\% &\\
\% &   \hspace{-15pt} General purpose sum of squares program (SOSP) solver:
\end{tabular}

\hspace{-20pt}\begin{tabular}{llcl}
\% &SOSPROGRAM &--- &Initialize a new SOSP.\\
\% &SOSDECVAR &--- &Declare new decision variables in an SOSP.\\
\% &SOSPOLYVAR &--- &Declare a new polynomial variable in an SOSP.\\
\% &SOSSOSVAR & --- &Declare a new sum of squares variable in an SOSP.\\
\% & SOSPOLYMATRIXVAR & --- & Declare a new matrix of polynomial variables in an SOSP.\\
\% & SOSSOSMATRIXVAR & --- & Declare a new matrix of sum of squares polynomial\\
\% & &  & variables in an SOSP.\\
\% &SOSPOSMATR       & --- & Declare a new positive semidefinite matrix variable in an SOSP\\
\% &SOSPOSMATRVAR    & --- & Declare a new symbolic positive semidefinite matrix \\
\% & & &                      variable in an SOSP\\
\% & SOSQUADVAR & --- & Declare a polynomial/SOS decision variable \\
\% & &  & with customized structure.\\
\% &SOSEQ & --- &Add a new equality constraint to an SOSP.\\
\% &SOSINEQ & --- &Add a new inequality constraint to an SOSP.\\
\% & SOSMATRIXINEQ & --- & Add a new matrix inequality constraint to an SOSP.\\
\% &SOSSETOBJ & --- &Set the objective function of an SOSP.\\
\% &SOSSOLVE & --- &Solve an SOSP.\\
\% &SOSGETSOL &--- &Get the solution from a solved SOSP.
\end{tabular}

\hspace{-20pt}\begin{tabular}{ll}
\% &\\
\% &   \hspace{-15pt} Customized functions:
\end{tabular}

\hspace{-20pt}\begin{tabular}{llcl}
\% &FINDSOS &--- &Find a sum of squares decomposition of a given polynomial\\
\% & &  &or a given polynomial matrix.\\
\% &FINDLYAP& --- &Find a Lyapunov function for a dynamical system.\\
\% &FINDBOUND& --- &Find a global/constrained lower bound for a polynomial.\\
\end{tabular}

\hspace{-20pt}\begin{tabular}{ll}
\% &\\
\% &   \hspace{-15pt} Demos:
\end{tabular}

\hspace{-20pt}\begin{tabular}{llcl}
\% &SOSDEMO1 and SOSDEMO1S &--- &Sum of squares test.\\
\% &SOSDEMO2 and SOSDEMO2S &--- &Lyapunov function search.\\
\% &SOSDEMO3 and SOSDEMO3S &--- &Bound on global extremum.\\
\% &SOSDEMO4 and SOSDEMO4S &--- &Matrix copositivity.\\
\% &SOSDEMO5 and SOSDEMO5S &--- &Upper bound for the structured singular value mu.\\
\% &SOSDEMO6 and SOSDEMO6S &--- &MAX CUT.\\
\% &SOSDEMO7S &--- &Chebyshev polynomials.\\
\% &SOSDEMO8S &--- &Bound in probability.\\
\% &SOSDEMO9 and SOSDEMO9S & --- & Sum of squares matrix decomposition.\\
\% &SOSDEMO10 and SOSDEMO10S & --- &Set containment.\\
\end{tabular}}

\verbatiminput{Contents.m}

\bibliographystyle{abbrv}
\bibliography{sostools}
\clearpage
\end{document}